\documentclass{amsart} 

\usepackage{amssymb,latexsym} 
 
\numberwithin{equation}{section} 
 
\newtheorem{theorem}{Theorem}[section] 
\newtheorem{proposition}[theorem]{Proposition} 
 
\newtheorem{corollary}[theorem]{Corollary} 
\newtheorem{remark}[theorem]{Remark} 
\newtheorem{lemma}[theorem]{Lemma} 
\newtheorem{example}[theorem]{Example} 
\newtheorem{definition}[theorem]{Definition}

\def\proof{\smallskip\noindent {\bf Proof. }} 
\def\endproof{\hfill$\square$\medskip}

\newcommand{\rem}[1]{\left\langle#1\right\rangle} 
 
\newcommand{\overunder}[2]{ 
\!\begin{array}{c} 
\scriptstyle{#1}\\[-.1in] 
-\!\!\!-\!\!\!-\\[-.1in] 
\scriptstyle{#2} 
\end{array} 
\! 
} 
\newcommand{\wideoverunder}[2]{ 
\!\begin{array}{c} 
\scriptstyle{#1}\\[-.1in] 
-\!\!\!-\!\!\!-\!\!\!-\!\!\!-\!\!\!-\!\!\!-\\[-.1in] 
\scriptstyle{#2} 
\end{array} 
\! 
} 
\newcommand{\overunderarrow}[2]{ 
\!\begin{array}{c} 
\scriptstyle{#1}\\[-.1in] 
-\!\!\!-\!\!\!\to\\[-.1in] 
\scriptstyle{#2} 
\end{array} 
\! 
} 
\newcommand{\wideoverunderarrow}[2]{ 
\!\begin{array}{c} 
\scriptstyle{#1}\\[-.1in] 
-\!\!\!-\!\!\!-\!\!\!-\!\!\!-\!\!\!-\!\!\!\to\\[-.1in] 
\scriptstyle{#2} 
\end{array} 
\! 
}

\def\AA{\mathcal{A}} 
\def\AAA{\mathbb{A}} 
 
\def\MM{\mathcal{M}} 
\def\PP{\mathbb{P}} 
\def\FFcal{\mathcal{F}} 
 
\def\TT{\mathbb{T}} 
\def\ZZ{\mathbb{Z}} 
\def\CC{\mathbb{C}} 
\def\RR{\mathbb{R}} 
 
\def\LL{\mathcal{L}}

\def\Trop{{\rm Trop}}

\newcommand{\bmat}[4]{\left[\!\!\begin{array}{cc} 
#1 & #2 \\ #3 & #4 \\ 
\end{array}\!\!\right]}

\begin{document} 
 
\title[Cluster algebras~I]{Cluster algebras~I: Foundations}

%    Information for first author~/www/Combin/winter01/ 
\author{Sergey Fomin} 
\address{Department of Mathematics, University of Michigan, 
Ann Arbor, MI 48109, USA} \email{fomin@math.lsa.umich.edu}

\author{Andrei Zelevinsky} 
\address{\noindent Department of Mathematics, Northeastern University, 
  Boston, MA 02115} 
\email{andrei@neu.edu} 
 
\begin{abstract}
In an attempt to create an algebraic framework for 
dual canonical bases and total positivity in semisimple groups,
we initiate the study of a new class of commutative algebras. 
\end{abstract}

\date{April 13, 2001} 
 
\dedicatory{To the memory of Sergei Kerov} 
 
\thanks{The authors were supported in part 
by NSF grants \#DMS-0049063, \#DMS-0070685 (S.F.), and  \#DMS-9971362 
(A.Z.). 
}

%    General info 
% \subjclass{Primary 14E05, % Rational and birational maps 
% Secondary 
% 12E05% Polynomials (irreducibility, etc.)
% 14E05% Rational and birational maps
% 14M15% Grassmannians, Schubert varieties, flag manifolds
% 14M17% Homogeneous spaces and generalizations 
% 17B37% Quantum groups (quantized enveloping algebras) and related
% deformations  
% 17B20% Simple, semisimple, reductive (super)algebras (roots)
% 17B67% Kac-Moody (super)algebras (structure and representation theory)
% 22E46% Semisimple Lie groups and their representations
% 22E65% Infinite-dimensional Lie groups and their Lie algebras
% 20G42% Quantum groups (quantized function algebras) and their representations
% 20G05% Representation theory [of linear algebraic groups] 
% } 
 
 
%{\ }

%\vspace{-.4in}
 
\maketitle 

\tableofcontents 
 
\section{Introduction} 
\label{sec:intro}

In this paper, we initiate the study of a new class of 
algebras, which we call \emph{cluster algebras}. 
Before giving precise definitions, we 
present some of the main features of these algebras. For any positive integer $n$, 
a cluster algebra $\AA$ of rank $n$ is a commutative ring with 
unit and no zero divisors, equipped with a distinguished family of generators 
called \emph{cluster variables}. The set of cluster variables is 
the (non-disjoint) union of a distinguished collection of 
$n$-subsets called \emph{clusters}. 
These clusters have the 
following \emph{exchange property}: for any cluster $\bf x$ and 
any element $x \in \bf x$, there is another cluster obtained from 
$\bf x$ by replacing $x$ with an element $x'$ related to $x$ by a 
\emph{binomial exchange relation} 
\begin{equation} 
\label{eq:binomial exchange general} x x' = M_1 + M_2 \ , 
\end{equation} 
where $M_1$ and $M_2$ are two monomials without common divisors in 
the $n-1$ variables $\textbf{x} - \{x\}$. Furthermore, any two 
clusters can be obtained from each other by a sequence of 
exchanges of this kind.

The prototypical example of a cluster algebra of rank 1 
is the coordinate ring $\AA = \CC [SL_2]$ of the group $SL_2$, viewed 
%as a cluster algebra 
in the following way. 
Writing a generic element of $SL_2$ as 
$\bmat{a}{b}{c}{d},$ 
we consider the entries $a$ and $d$ as cluster variables, 
and the entries $b$ and $c$ as scalars. 
There are just two clusters $\{a\}$ and $\{d\}$, 
and $\AA$ is the algebra over the polynomial ring 
$\CC [b,c]$ generated by the cluster variables $a$ and $d$ 
subject to the binomial exchange relation 
$$ad = 1 + b c \ .$$ 
 
Another important incarnation of a cluster algebra of rank 1 
is the coordinate ring $\AA = \CC [SL_3/N]$ of the base affine space 
of the special linear group $SL_3$; here $N$ is the maximal unipotent subgroup 
of $SL_3$ consisting of all unipotent upper triangular matrices. 
Using the standard notation $(x_1, x_2, x_3, x_{12}, x_{13}, 
x_{23})$ for the Pl\"ucker coordinates on $SL_3/N$, we view $x_2$ and $x_{13}$ 
as cluster variables; 
then $\AA$ is the algebra over the polynomial ring 
$\CC [x_1,x_3, x_{12}, x_{13}]$ generated by the two cluster variables 
$x_2$ and $x_{13}$ 
subject to the binomial exchange relation 
$$x_2 x_{13} = x_1 x_{23} + x_3 x_{12} \ .$$ 
This form of representing the algebra $\CC [SL_3/N]$ is closely related to the 
choice of a linear basis in it consisting of all monomials in 
the six Pl\"ucker coordinates which are not divisible by $x_2 x_{13}$. 
This basis was introduced and studied in \cite{gz86} under the name 
``canonical basis." 
As a representation of $SL_3$, the space $\CC [SL_3/N]$ 
is the multiplicity-free direct sum of all irreducible 
finite-dimensional representations, and each of the components 
is spanned by a part of the above basis. 
Thus, this construction provides a ``canonical" basis 
in every irreducible finite-dimensional representation of~$SL_3$. 
After Lusztig's work \cite{lu1}, this basis had been 
recognized as (the classical limit at $q \to 1$ 
of) the \emph{dual canonical basis}, i.e., the basis 
in the $q$-deformed algebra $\CC_q [SL_3/N]$ which is dual to 
Lusztig's canonical basis in the appropriate $q$-deformed 
universal enveloping algebra (a.k.a. quantum group). 
The dual canonical basis in the space $\CC [G/N]$ was later 
constructed explicitly for a few 
other classical groups $G$ of small rank: for $G = Sp_4$ in 
\cite{rz} and for $G = SL_4$ in \cite{bz93}. 
In both cases, $\CC [G/N]$ can be seen to be a cluster algebra: 
there are 6 clusters of size 2 for $G = Sp_4$, and 14 clusters of size 
3 for for $G = SL_4$. 
 
We conjecture that the above examples can be extensively 
generalized: for any simply-connected connected semisimple group 
$G$, the coordinate rings $\CC [G]$ and $\CC [G/N]$, as well as 
coordinate rings of many other interesting varieties related  
to~$G$, have a natural structure of a cluster algebra. 
This structure should serve as an algebraic framework for the study 
of ``dual canonical bases" in these coordinate rings and 
their $q$-deformations. 
In particular, we conjecture that all monomials in the variables 
of any given cluster (the \emph{cluster monomials}) 
belong to this dual canonical basis. 
 
A particularly nice and well understood example of a cluster 
algebra of an arbitrary rank $n$ is the homogeneous coordinate ring 
$\CC [Gr_{2,n+3}]$ of the Grassmannian of $2$-dimensional 
subspaces in~$\CC^{n+3}$. 
This ring is generated by the Pl\"ucker coordinates $[ij]$, 
for $1 \leq i < j \leq n+3$, subject to the relations 
$$[ik][jl] = [ij] [kl] + [il] [jk] \ ,$$ 
for all $i < j < k < l$. 
It is convenient to identify the indices $1, \dots, n+3$ with the 
vertices of a convex $(n+3)$-gon, and  
the Pl\"ucker coordinates with its sides and diagonals. 
We view the sides $[12],[23], \dots, [n+2, n+3], [1,n+3]$ as 
scalars, and the diagonals as cluster variables. 
The clusters are the maximal families of pairwise non-crossing 
diagonals; thus, they are in a natural bijection with the 
triangulations of this polygon. 
It is known that the cluster monomials form a linear basis 
in~$\CC [Gr_{2,n+3}]$. 
To be more specific, we note that this ring is naturally identified with the ring of 
polynomial $SL_2$-invariants of an $(n+3)$-tuple of points in~$\CC^2$. 
Under this isomorphism, the basis of cluster monomials corresponds to the basis considered 
in \cite{kungrota84, stur}. (We are grateful to Bernd Sturmfels for bringing 
these references to our attention.)  
%This basis contains all cluster monomials, 
%Again, it can be seen from this description that $\CC [Gr_{2,n+3}]$ 
%is a cluster algebra of rank $n$: 
%and the clusters are in a 
%natural bijection with triangulations of a convex $(n+3)$-gon. 
%This example will be placed in a more general context in the sequel 
%to this paper. 
 
An essential feature of the exchange relations 
(\ref{eq:binomial exchange general}) is that the right-hand side does not involve subtraction. 
Recursively applying these relations, one can represent any cluster variable as a 
\emph{subtraction-free} rational expression in the 
variables of any given cluster. 
This positivity property is consistent with 
a remarkable connection between canonical bases and the 
theory of total positivity, discovered by 
G.~Lusztig \cite{lusztig-quantum,lusztig-reductive}. 
Generalizing the classical concept of totally positive matrices, 
he defined totally positive elements in any reductive 
group $G$, and proved that all elements of the dual canonical 
basis in $\CC[G]$ take positive values at them. 
 
It was realized in \cite{lusztig-reductive,fz-jams} 
that the natural geometric framework 
for total positivity is given by \emph{double Bruhat cells},  
the intersections of cells of the Bruhat decompositions with 
respect to two opposite Borel subgroups. 
Different aspects of total positivity in double Bruhat cells were explored 
by the authors of the present paper and their collaborators in 
\cite{bfz96,bz97,bz01,fz-jams,fz-intelligencer, 
fz-proceedings,lekzel98,ssvz,z-imrn}. 
The binomial exchange relations of the form 
(\ref{eq:binomial exchange general}) played a crucial role 
in these studies. 
It was the desire to explain the ubiquity of these relations 
and to place them in a proper context that led us to 
the concept of cluster algebras. 
The crucial step in this direction was made 
in \cite{z-imrn}, where a family of clusters and exchange relations 
was explicitly constructed in the coordinate ring of an arbitrary double Bruhat cell. 
However, this family was not \emph{complete}: 
in general, some clusters were missing, 
%not any 
%two of them could be joined by a sequence of exchanges, 
and not any member of a cluster could be exchanged from it. 
Thus, we started looking for a natural way to 
``propagate" exchange relations from one cluster to another. 
The concept of cluster algebras is the result of this investigation. 
We conjecture that the coordinate ring of any double Bruhat cell is a cluster algebra.  
 
This article, in which we develop the foundations of 
the theory, is conceived as the first in a forthcoming series. 
%Here we introduce the basic notions of the theory
%and develop its most fundamental concepts.
%and obtain some of their algebraic consequences. 
We attempt to make the exposition elementary and 
self-contained; in particular, no knowledge of semisimple groups, quantum groups or 
total positivity is assumed on the part of the reader. 
 
One of the main structural features of cluster algebras established in the present paper 
is the following  
%The main result of this paper is that cluster algebras possess 
%an important algebraic feature that we call the 
\emph{Laurent phenomenon}: any cluster variable $x$ 
viewed as a  rational function in the variables of any given cluster 
is in fact a Laurent polynomial. 
This property is quite surprising: in most cases, 
the numerators of these Laurent polynomials 
contain a huge number of monomials, and the numerator for $x$ 
moves into the denominator when we compute the cluster variable $x'$ obtained 
from $x$ by an exchange (\ref{eq:binomial exchange general}). 
The magic of the Laurent phenomenon is that, at every stage of 
this recursive process, a cancellation will inevitably occur, 
leaving a single monomial in the denominator. 
 
In view of the positivity property discussed above, it is natural to expect that all Laurent polynomials 
for cluster variables will have \emph{positive} coefficients. 
This seems to be a rather deep property; our present methods do not 
provide a proof of it.  
 
On the bright side, by a modification of the method developed here, 
it is possible to establish the Laurent phenomenon in many different 
situations spreading beyond the cluster algebra framework. 
We explore these situations in a separate paper \cite{fz-Laurent}. 
 
The paper is organized as follows. 
Section~\ref{sec:setup} contains an axiomatic definition, first examples and 
the first structural properties of cluster algebras. 
One of the technical difficulties in setting up the foundations involves 
the concept of an \emph{exchange graph} 
whose vertices correspond to clusters, and the edges to exchanges among them. 
It is convenient to begin by taking the $n$-regular tree $\TT_n$ as our underlying graph. 
This tree can be viewed as a universal cover for the actual exchange graph, whose appearance is postponed 
until Section~\ref{sec:cluster graph}. 
 
The Laurent phenomenon is established in Section~\ref{sec:laurent}. 
In Sections~\ref{sec:exponents} and \ref{sec:coefficients}, 
we scrutinize the main definition, 
obtain useful reformulations, and introduce some important 
classes of cluster algebras.

Section~\ref{sec:rank 2} contains a detailed analysis of cluster algebras of rank 2. 
This analysis exhibits deep and somewhat mysterious connections 
between cluster algebras and Kac-Moody algebras. 
This is just the tip of an iceberg: these connections 
will be further explored (for cluster algebras of an arbitrary rank) in the sequel to this paper. 
The main result of this sequel is a complete classification 
of cluster algebras of \emph{finite type}, i.e., those with finitely many 
distinct clusters; cf.\ Example~\ref{ex:rank 2 exchange graph}. 
This classification turns out to be yet another instance of the famous Cartan-Killing 
classification.

\medskip 
 
\noindent \textsc{Acknowledgments.} This work began in May 2000 
when the authors got together as participants in 
the program ``Representation Theory-2000," organized by Victor Kac 
and Alexander A.~Kirillov at the Erwin Schr\"odinger International Institute 
for Mathematical Physics in Vienna, Austria. We thank the organizers 
for inviting us, and we are grateful to Peter Michor 
and the staff of the Institute for creating ideal working conditions. 
The paper was finished at the Isaac Newton Institute for Mathematical Sciences 
in Cambridge, UK, whose support is gratefully acknowledged.

It was during our stay in Vienna that we learned of the terminal illness 
of Sergei Kerov, an outstanding mathematician and a good friend of ours. 
He passed away on July 30, 2000. We dedicate this paper to his memory.

%\section*{Basic notions} 
 
\section{Main definitions} 
\label{sec:setup}

Let $I$ be a finite set of size $n$; the standard choice 
will be $I = [n]=\{1,2,\dots,n\}$. 
Let $\TT_n$ denote the \emph{$n$-regular tree,} whose edges are 
labeled by the elements of $I$, 
so that the $n$ edges emanating from each vertex receive 
different labels. 
By a common abuse of notation, we will sometimes denote by $\TT_n$ 
the set of the tree's vertices. 
We will write $t \overunder{i}{} t'$ if vertices 
$t,t'\in\TT_n$ are joined by an edge labeled by~$i$.

To each vertex $t\in \TT_n$, we will associate a \emph{cluster} of $n$ 
generators (``variables") ${\bf x}(t)=(x_i(t))_{i \in I}$. 
All these variables will commute with each other and satisfy the following 
\emph{exchange relations,} for 
every edge $t \overunder{j}{} t'$ in $\TT_n$:
\begin{eqnarray} 
\label{eq:exchange1} 
&&\text{$x_i(t)=x_i(t')$\quad for any $i\neq j$;}\\[.1in] 
\label{eq:exchange2} 
&&\text{$x_j(t)\,x_j(t')= M_j (t)({\bf x}(t)) + 
M_j (t')({\bf x}(t'))$.} 
\end{eqnarray} 
Here $M_j (t)$ and $M_j (t')$ are two monomials in 
the $n$ variables $x_i \, , i \in I$; we think of these monomials 
as being associated with the two ends of the edge $t \overunder{j}{} t'$. 
 
To be more precise, let $\PP$ be an abelian group without torsion, written 
multiplicatively. 
We call $\PP$ the \emph{coefficient group}; 
%The group ring $\ZZ\PP$ has no zero divisors; 
%let $\FF$ denote its quotient field. 
a prototypical example is a free abelian group of finite rank. 
Every monomial $M_j (t)$ in (\ref{eq:exchange2}) will have the form 
\begin{equation} 
\label{eq:M(t)} 
M_j(t)=p_j(t)\, \prod_{i \in I} x_i^{b_i} \, , 
\end{equation} 
for some coefficient $p_j(t)\in\PP$ and some nonnegative integer 
exponents $b_i$. 
%Informally, we think of $M_j (t)$ as being associated with the $t$-end 
%of the edge $t \overunder{j}{} t'$. 
 
The monomials $M_j(t)$ must satisfy certain conditions (axioms). 
To state them, we will need a little preparation. 
Let us write $P\mid Q$ to denote that a polynomial $P$ divides a 
polynomial~$Q$. 
Accordingly, $x_i\!\mid\! M_j(t)$ means that the monomial $M_j(t)$ 
contains the variable~$x_i$. 
For a rational function $F=F(x,y,\dots)$, 
the notation $F|_{x\leftarrow g(x,y,\dots)}$ 
will denote the result of substituting $g(x,y,\dots)$ for $x$ into~$F$. 
To illustrate, if $F(x,y)=xy$, then $F|_{x\leftarrow \frac{y}{x}}=\frac{y^2}{x}$. 
 
\begin{definition} 
\label{def:exchange pattern} 
{\rm An \emph{exchange pattern} on $\TT_n$ with coefficients in 
$\PP$ is a family of monomials $\MM=(M_j(t))_{t\in\TT_n, \ j \in I}$ of the form 
(\ref{eq:M(t)}) satisfying the following four 
axioms: 
\begin{eqnarray} 
\label{eq:CA1} 
&&\text{If $t \in\TT_n$, then $x_j \nmid M_j(t)$.} 
\\ 
\label{eq:CA2} 
&&\text{If $t_1 \overunder{j}{} t_2$ and $x_i\mid M_j(t_1)$, then 
$x_i\nmid M_j(t_2)$.} \\ 
\label{eq:CA3} 
&&\text{If $t_1 \overunder{i}{} t_2 \overunder{j}{} t_3\,$, then 
$~~x_j\mid M_i(t_1)~~$ if and only if  $~~x_i\mid M_j(t_2)~~$.}\\ 
\label{eq:CA4} 
&&\text{Let $~t_1 \overunder{i}{} t_2 \overunder{j}{} t_3 
  \overunder{i}{} t_4~$. Then 
$~~\displaystyle\frac{M_i(t_3)}{M_i(t_4)} 
= \displaystyle\frac{M_i(t_2)}{M_i(t_1)}\Bigl|_{x_j\leftarrow 
{M_0}/{x_j}}\,$,} 
\\ 
\nonumber 
& &\text{where $M_0=(M_j(t_2)+M_j(t_3))|_{x_i=0}\,$.} 
\end{eqnarray} 
} 
\end{definition} 
 
We note that in the last axiom, the substitution $x_j\leftarrow 
\frac{M_0}{x_j}$ is effectively monomial, since 
in the event that neither $M_j(t_2)$ nor $M_j(t_3)$ 
contain $x_i$, condition (\ref{eq:CA3}) requires that both $M_i(t_2)$ 
and $M_i(t_1)$ do not depend on $x_j$, thus making the whole 
substitution irrelevant. 
%redundant. 
 
One easily checks that axiom (\ref{eq:CA4}) is invariant under the 
``flip'' $t_1\leftrightarrow t_4$, $t_2\leftrightarrow t_3$, 
so no restrictions are added if we apply it ``backwards.'' 
The axioms also imply at once that setting 
\begin{equation} 
\label{eq:exchange pattern symmetry} 
M'_j (t) = M_j (t') 
\end{equation} 
for every edge $t \overunder{j}{} t'$, we obtain another exchange 
pattern $\MM'$; this gives a natural involution 
$\MM \to \MM'$ on the set of all exchange patterns. 
 
\begin{remark} 
\label{rem:degree of freedom} 
{\rm Informally speaking, axiom (\ref{eq:CA4}) describes 
the propagation of an exchange pattern along the edges of $\TT_n$. 
More precisely, let us fix the $2n$ exchange monomials for all edges emanating from 
a given vertex $t$. 
This choice uniquely determines the ratio ${M_i(t')}/{M_i(t'')}$ 
for any vertex $t'$ adjacent to $t$ and any edge $t' \overunder{i}{} t''$ 
(to see this, take $t_2 = t$ 
and $t_3 = t'$ in (\ref{eq:CA4}), and allow $i$ to vary). 
In view of (\ref{eq:CA2}), this ratio in turn uniquely determines 
the exponents of all variables $x_k$ in both monomials 
$M_i(t')$ and $M_i(t'')$. 
There remains, however, one degree of freedom in determining 
the coefficients $p_i(t')$ and $p_i(t'')$ because only their ratio 
is prescribed by (\ref{eq:CA4}). 
In Section~\ref{sec:coefficients} we shall introduce an 
important class of \emph{normalized} exchange patterns for which 
this degree of freedom disappears, and so the whole pattern  
is uniquely determined by the $2n$ monomials associated with edges 
emanating from a given vertex.} 
\end{remark} 

Let $\ZZ \PP$ denote the group ring of $\PP$ with integer coefficients. 
For an edge $t \overunder{k}{} t'$, we refer to the binomial 
$P=M_k(t)+M_k(t') \in \ZZ \PP [x_i : i \in I]$ as the 
\emph{exchange  polynomial}. 
We will write $~t \overunder{}{P} t'~$ or $~t \overunder{k}{P} 
t'~$ to indicate this fact. 
Note that, in view of the axiom (\ref{eq:CA1}), the right-hand 
side of the exchange relation (\ref{eq:exchange2}) can be written 
as $P({\bf x}(t))$, which is the same as $P({\bf x}(t'))$. 
%Since $\PP$ is torsion-free, $-1\notin\PP$ and therefore $P\neq 0$. 

Let $\MM$ be an exchange pattern on $\TT_n$ with coefficients in $\PP$. 
Note that since $\PP$ is torsion-free, the ring $\ZZ \PP$ 
has no zero divisors. 
%let $\FF$ denote its quotient field. 
For every vertex $t\in \TT_n$, let $\FFcal(t)$ denote the field of rational functions 
in the cluster variables $x_i (t)$, $i \in I$, with coefficients in~$\ZZ \PP$. 
For every edge $t \overunder{j}{P} t'$, we define a 
$\ZZ \PP$-linear field isomorphism 
%birational \emph{transition map} 
$R_{tt'}:\FFcal(t')\to\FFcal(t)$ by 
\begin{equation} 
\label{eq:transition1} 
\begin{array}{l} 
R_{tt'}(x_i(t'))=x_i(t)\quad \text{for $i\neq k$;}\\[.1in] 
R_{tt'}(x_k(t'))=\displaystyle\frac{P({\bf x}(t))}{x_k(t)}\,. 
\end{array} 
\end{equation} 
Note that property (\ref{eq:CA1}) ensures that 
$R_{t't}=R_{tt'}^{-1}$. 
%We then define, for any path 
%$~ 
%t_0 \overunder{}{} t_1 \overunder{}{} \cdots 
%  \overunder{}{} t_k \,, 
%$ 
%the transition map 
%\begin{equation} 
%\label{eq:transition2} 
%R_{t_0 t_k} = 
%R_{t_0 t_{1}}\circ \cdots\circ R_{t_{k-1} t_k} 
%\,:\,\FFcal(t_k)\to\FFcal(t_0) 
%\,. 
%\end{equation} 
The \emph{transition maps} $R_{tt'}$ enable us to identify all the 
fields $\FFcal(t)$ with each other. 
We can then view them as a single field $\FFcal$ that 
contains all the elements $x_i(t)$, for all $t\in\TT_n$ and $i\in I$. 
Inside~$\FFcal$, these elements satisfy the exchange relations 
(\ref{eq:exchange1})--(\ref{eq:exchange2}).

\begin{definition} 
\label{def:cluster-algebra} 
{\rm 
Let $\AAA$ be a subring with unit in $\ZZ \PP$ containing all coefficients 
$p_i (t)$ for $i \in I$ and $t\in\TT_n$. 
The \emph{cluster algebra} $\AA=\AA_{\AAA} (\MM)$ of rank $n$ over 
$\AAA$ associated with an exchange pattern $\MM$ is the %commutative 
$\AAA$-subalgebra with unit in  
$\FFcal$ generated by the union of all clusters~${\bf x}(t)$, 
for $t\in\TT_n$. 
} 
\end{definition} 
 
The smallest possible ground ring $\AAA$ is the subring of 
$\ZZ \PP$ generated by all the coefficients $p_i (t)$; 
the largest one is $\ZZ \PP$ itself. 
%A natural intermediate choice is to take $\AAA = \ZZ \PP$. 
An intermediate choice of $\AAA$ appears in Proposition~\ref{pr:reduction} below.

Since $\AA$ is a subring of a field $\FFcal$, it is a commutative ring with no zero divisors. 
We also note that if $\MM'$ is obtained from $\MM$ by the involution 
(\ref{eq:exchange pattern symmetry}), then the cluster algebra $\AA_\AAA (\MM')$ 
is naturally identified with $\AA_\AAA (\MM)$.

\begin{example} 
\label{ex:rank 1} 
{\rm Let $n = 1$. 
The tree $\TT_1$ has only one edge $t \overunder{1}{} t'$. 
The corresponding cluster algebra $\AA$ has two generators 
$x = x_1 (t)$ and $x' = x_1 (t')$ satisfying the exchange relation 
$$x x' = p + p' \ ,$$ 
where $p$ and $p'$ are arbitrary elements of the coefficient 
group $\PP$. 
In the ``universal" setting, we take $\PP$ to be the free abelian 
group generated by $p$ and $p'$. 
Then the two natural choices for the ground ring $\AAA$ are 
the polynomial ring $\ZZ [p,p']$, and the Laurent polynomial ring 
$\ZZ \PP = \ZZ [p^{\pm 1},{p'}^{\ \pm 1}]$. 
All other realizations of $\AA$ can be viewed as 
specializations of the universal one. 
Despite the seeming triviality of this example, it covers several 
important algebras: the coordinate ring of each of the varieties 
$SL_2$, $Gr_{2,4}$ and $SL_3/B$ (cf. Section~\ref{sec:intro}) 
is a cluster algebra of rank~$1$, for an appropriate choice of 
$\PP$, $p$, $p'$ and $\AAA$. 
} 
\end{example} 
 
\begin{example} 
\label{ex:rank 2} 
{\rm Consider the case $n = 2$. 
The tree $\TT_2$ is shown below: 
\begin{equation} 
\label{eq:tree-2} 
\cdots \overunder{1}{} 
%\bullet\wideoverunder{2}{a_{-1}x_1+b_{-1}} 
%t_0 \wideoverunder{1}{} 
t_0 \wideoverunder{2}{} 
t_1 \wideoverunder{1}{} 
t_2 \wideoverunder{2}{} 
t_3 \wideoverunder{1}{} 
t_4 \overunder{2}{} 
\cdots \,. 
\end{equation} 
Let us denote the cluster variables as follows: 
\[ 
y_1=x_1(t_0)=x_1(t_1), \qquad 
y_2=x_2(t_1)=x_2(t_2),  \qquad 
y_3=x_1(t_2)=x_1(t_3), \dots 
\] 
(the above equalities among the cluster variables follow from 
(\ref{eq:exchange1})). 
Then the clusters look like 
\[ 
\cdots \overunder{1}{}\hspace{-.1in} 
\begin{array}{c}\scriptstyle y_1,y_0\\ \bullet\\ t_0\end{array} 
 \hspace{-.1in}\wideoverunder{2}{}\hspace{-.1in} 
\begin{array}{c}\scriptstyle y_2,y_1\\ \bullet\\ t_1\end{array} 
 \hspace{-.1in}\wideoverunder{1}{}\hspace{-.1in} 
\begin{array}{c}\scriptstyle y_3,y_2\\ \bullet\\ t_2\end{array} 
 \hspace{-.1in}\wideoverunder{2}{}\hspace{-.1in} 
\begin{array}{c}\scriptstyle y_4,y_3\\ \bullet\\ t_3\end{array} 
 \hspace{-.1in}\wideoverunder{1}{}\hspace{-.1in} 
\begin{array}{c}\scriptstyle y_5,y_4\\ \bullet\\ t_4\end{array} 
% \hspace{-.1in}\wideoverunder{2}{}\hspace{-.1in} 
%\begin{array}{c}\scriptstyle y_3,y_4\\ \bullet\\ t_4\end{array} 
% \hspace{-.1in}\wideoverunder{1}{}\hspace{-.1in} 
%\begin{array}{c}\scriptstyle y_5,y_4\\ \bullet\\ t_5\end{array} 
 \hspace{-.1in}\overunder{2}{} 
\cdots \,. 
\] 
We claim that the exchange relations (\ref{eq:exchange2}) can be 
written in the following form: 
\begin{equation} 
\label{eq:exchange rank 2} 
\begin{array}{ll} 
y_{0} y_2 = q_1 y_1^{b}+r_1, \quad & y_{1} y_3 = q_2 y_2^{c}+r_2, \\[.1in] 
y_{2} y_4 = q_3 y_3^{b}+r_3, & y_{3} y_5 = q_4 y_4^{c}+r_4, \dots , 
\end{array} 
\end{equation} 
where the integers $b$ and $c$ are either both positive 
or both equal to~$0$, and the coefficients $q_m$ and $r_m$ are 
elements of $\PP$ satisfying the relations 
\begin{equation} 
\label{eq:qqr=rr} 
\begin{array}{ll} 
q_0 q_2 r_1^c = r_0 r_2, \quad & q_1 q_3 r_2^b = r_1 r_3, \\[.1in] 
q_2 q_4 r_3^c = r_2 r_4, & q_3 q_5 r_4^b = r_3 r_5, \, \dots . 
\end{array} 
\end{equation} 
Furthermore, any such choice of parameters $b$, $c$, $(q_m)$, $(r_m)$ 
results in a well defined cluster algebra of rank~2. 
 
To prove this, we notice that, in view of 
(\ref{eq:CA1})--(\ref{eq:CA2}), both monomials $M_2 (t_0)$ and $M_2 
(t_1)$ do not 
contain the variable $x_2$, and at most one of them contains $x_1$. 
If $x_1$ enters neither $M_2 (t_0)$ nor $M_2 (t_1)$, then 
these two are simply elements of~$\PP$. 
But then  (\ref{eq:CA3}) forces \emph{all} 
monomials~$M_i (t_m)$ to be elements of $\PP$, while (\ref{eq:CA4}) implies 
that it is possible to give the names $q_m$ and $r_m$ to the two monomials corresponding to 
each edge $t_m \overunder{}{} t_{m+1}$ so that 
(\ref{eq:exchange rank 2})--(\ref{eq:qqr=rr}) hold with $b = c = 0$. 
 
Next, consider the case when precisely one of the 
monomials $M_2 (t_0)$ and $M_2 (t_1)$ contains $x_1$. 
Applying if necessary the involution (\ref{eq:exchange pattern symmetry}) to 
our exchange pattern, we may assume that $M_2 (t_0) = q_1 x_1^b$ 
and $M_2 (t_1) = r_1$ for some positive integer $b$ and some $q_1, r_1 \in \PP$. 
Thus, the exchange relation associated to the edge $t_0 
\overunder{2}{} t_1$ takes the form $y_{0} y_2 = q_1 y_1^{b}+r_1$. 
By (\ref{eq:CA3}), we have $M_1 (t_1) = q_2 x_2^c$ 
and $M_1 (t_2) = r_2$ for some positive integer $c$ and some $q_2, r_2 \in \PP$. 
Then the exchange relation for the edge $t_1 
\overunder{1}{} t_2$ takes the form $y_{1} y_3 = q_2 
y_2^{c}+r_2$. 
At this point, we invoke (\ref{eq:CA4}): 
$$ 
\frac{M_2(t_2)}{M_2(t_3)} 
= \frac{M_2(t_1)}{M_2(t_0)}\Bigl|_{x_1 \leftarrow r_2/x_1} 
= \frac{r_1}{q_1 
x_1^b}\Bigl|_{x_1 \leftarrow r_2/x_1} 
= \frac{r_1 x_1^b}{q_1 r_2^b} \ .$$ 
By (\ref{eq:CA2}), we have $M_2 (t_2) = q_3 x_1^b$ and $M_2 (t_3) = r_3$ 
for some $q_3, r_3 \in \PP$ satisfying $q_1 q_3 r_2^b = r_1 r_3$. 
Continuing in the same way, we obtain all relations 
(\ref{eq:exchange rank 2})--(\ref{eq:qqr=rr}). 
 
For fixed $b$ and $c$, the ``universal" coefficient group $\PP$ is 
the multiplicative abelian group generated by the elements $q_m$ and $r_m$ for 
all $m \in \ZZ$ subject to the defining relations (\ref{eq:qqr=rr}). 
It is easy to see that this is a free abelian group of infinite rank. 
As a set of its free generators, one can choose any subset 
of $\{q_m, r_m: m \in \ZZ\}$ that contains  four generators 
$q_0, r_0, q_1, r_1$ and precisely one generator from each pair 
$\{q_m, r_m\}$ for $m \neq 0, 1$. 
 
A nice specialization of this setup is provided by the homogeneous
coordinate ring of the Grassmannian $Gr_{2,5}$. 
Recall (cf. Section~\ref{sec:intro}) that this ring is generated by the Pl\"ucker 
coordinates $[k,l]$, where $k$ and $l$ are distinct elements of 
the cyclic group $\ZZ / 5 \ZZ$. We shall write 
$\overline m = m \mod  5 \in \ZZ / 5 \ZZ$ 
for $m \in \ZZ$, and adopt the convention $[k,l] = [l,k]$; see 
Figure~\ref{fig:pentagon}. 
 
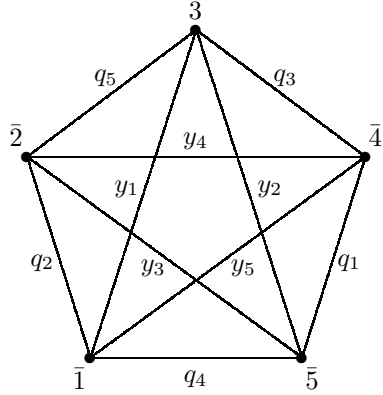
\begin{figure}[ht] 
\begin{center} 
\setlength{\unitlength}{4pt} 
\begin{picture}(40,34)(0,-2) 
%\thicklines 
  \put(6,0){\line(1,0){20}} 
  \put(0,19){\line(1,0){32}} 
  \qbezier(6,0)(11,15.5)(16,31) 
  \qbezier(26,0)(21,15.5)(16,31) 
  \qbezier(6,0)(3,9.5)(0,19) 
  \qbezier(26,0)(29,9.5)(32,19) 
  \qbezier(0,19)(13,9.5)(26,0) 
  \qbezier(32,19)(19,9.5)(6,0) 
  \qbezier(0,19)(8,25)(16,31) 
  \qbezier(32,19)(24,25)(16,31) 
 
  \put(6,0){\circle*{1}} 
  \put(26,0){\circle*{1}} 
  \put(0,19){\circle*{1}} 
  \put(32,19){\circle*{1}} 
  \put(16,31){\circle*{1}} 
 
\put(5,-2){\makebox(0,0){$\bar 1$}} 
\put(27,-2){\makebox(0,0){$\bar 5$}} 
\put(-1,21){\makebox(0,0){$\bar 2$}} 
\put(33,21){\makebox(0,0){$\bar 4$}} 
\put(16,33){\makebox(0,0){$\bar 3$}} 
\put(16,-2){\makebox(0,0){$q_4$}} 
\put(1.5,9){\makebox(0,0){$q_2$}} 
\put(30.5,9){\makebox(0,0){$q_1$}} 
\put(7.5,26.5){\makebox(0,0){$q_5$}} 
\put(24.5,26.5){\makebox(0,0){$q_3$}} 
\put(16,20.5){\makebox(0,0){$y_4$}} 
\put(23,16){\makebox(0,0){$y_2$}} 
\put(9.5,16){\makebox(0,0){$y_1$}} 
\put(12,8.5){\makebox(0,0){$y_3$}} 
\put(20.5,8.5){\makebox(0,0){$y_5$}} 
 
\end{picture} 
\end{center} 
\caption{The Grassmannian $Gr_{2,5}$} 
\label{fig:pentagon} 
\end{figure} 
 
The ideal of relations among the Pl\"ucker coordinates is generated by the 
relations 
$$[\overline m, \overline {m+2}] [\overline {m+1}, \overline {m+3}] = 
[\overline m, \overline {m+1}] [\overline {m+2}, \overline {m+3}] + 
[\overline m, \overline {m+3}] [\overline {m+1}, \overline {m+2}]$$ 
for $m \in \ZZ$. 
Direct check shows that these relations are a specialization 
of the relations (\ref{eq:exchange rank 2}), if we set 
$b = c = 1$, $y_m =  [\overline {2m-1}, \overline {2m+1}]$, 
$q_m = [\overline {2m-2}, \overline {2m+2}]$, and 
$r_m = [\overline {2m-2}, \overline {2m-1}][\overline {2m+1}, \overline {2m+2}] 
= q_{m-2} q_{m+2}$ for 
all $m \in \ZZ$. 
The coefficient group $\PP$ is the multiplicative free abelian group 
with $5$ generators $q_m$. 
It is also immediate that the elements $q_m$ and $r_m$ defined in this way satisfy 
the relations (\ref{eq:qqr=rr}).  
} 
\end{example} 
 
We conclude this section by introducing two important operations 
on exchange patterns: \emph{restriction} and \emph{direct product}. 
Let us start with restriction.  
Let $\MM$ be an exchange pattern of rank $n$ with an index set 
$I$ and coefficient group $\PP$. 
Let $J$ be a subset of size $m$ in $I$. 
Let us remove from $\TT_n$ all edges labeled by 
indices in $I - J$, and choose any connected component of 
the resulting graph. 
This component is naturally identified with $\TT_m$. 
Let $\MM'$ denote the restriction of $\MM$ to $\TT_m$, i.e., 
the collection of monomials $M_j (t)$ for all $j \in J$ and 
$t \in \TT_m$. 
Then $\MM'$ is an exchange pattern on $\TT_m$ whose 
coefficient group $\PP'$ is the direct product of $\PP$ with 
the multiplicative free abelian group with generators 
$x_i$, $i \in I - J$. 
We shall say that $\MM'$ is obtained from $\MM$ by restriction 
from $I$ to $J$.

\begin{proposition} 
\label{pr:reduction} 
Let $\AA=\AA_{\AAA} (\MM)$ be a cluster algebra of rank $n$. 
The $\AAA$-subalgebra of $\AA$ generated by %the union of clusters 
$\cup_{t \in \TT_m} {\bf x} (t)$ %for $t \in \TT_m$ 
is naturally identified with the cluster algebra $\AA_{\AAA'} (\MM')$, where 
$\AAA'$ is the polynomial ring $\AAA [x_i: i \in I - J]$. 
\end{proposition} 
 
\proof 
If $i \in I - J$ then (\ref{eq:exchange1}) implies 
that $x_i (t)$ stays constant as $t$ varies over $\TT_m$. 
Therefore, we can identify this variable 
with the corresponding generator $x_i$ of the coefficient group 
$\PP'$, and the statement follows. 
\endproof 
 
Let us now consider two exchange patterns $\MM_1$ and $\MM_2$ 
of ranks $n_1$ and $n_2$, respectively, 
with index sets $I_1$ and $I_2$, and coefficient groups $\PP_1$ and $\PP_2$. 
We will construct the exchange pattern $\MM = \MM_1 \times \MM_2$ 
(the direct product of $\MM_1$ and $\MM_2$) of rank $n = n_1 + n_2$, 
with the index set $I = I_1 \bigsqcup I_2$, and coefficient group 
$\PP = \PP_1 \times \PP_2$. 
Consider the tree $\TT_n$ whose edges are colored by $I$, and, 
for $\nu \in \{1,2\}$, let $\pi_\nu : \TT_n \to \TT_{n_\nu}$ 
be a map with the following property: if 
$t \overunder{i}{} t'$ in $\TT_n$, and $i \in I_\nu$ 
(resp., $i \in I - I_\nu = I_{3 - \nu}$), then 
$\pi_\nu (t) \overunder{i}{} \pi_\nu (t')$ in 
$\TT_{n_\nu}$ (resp., $\pi_\nu (t) = \pi_\nu (t')$). 
Clearly, such a map $\pi_\nu$ exists and is essentially unique: 
it is determined by specifying the image of any vertex of $\TT_n$. 
We now introduce the exchange pattern $\MM$ on $\TT_n$ by setting, 
for every $t \in \TT_n$ and $i \in I_\nu \subset I$, the monomial 
$M_i (t)$ to be equal to $M_i (\pi_\nu (t))$, the latter monomial 
coming from the exchange pattern $\MM_\nu$. 
The axioms (\ref{eq:CA1})--(\ref{eq:CA4}) for $\MM$ are checked 
directly.

\begin{proposition} 
\label{pr:tensor product} 
Let $\AA_1=\AA_{\AAA_1} (\MM_1)$ and 
$\AA_2=\AA_{\AAA_2} (\MM_2)$ be cluster algebras. 
Let $\MM = \MM_1 \times \MM_2$ and $\AAA = \AAA_1 \otimes \AAA_2$. 
Then the cluster algebra $\AA_{\AAA} (\MM)$ is canonically isomorphic 
to the tensor product of algebras $\AA_1 \otimes \AA_2$ 
(all tensor products are taken over $\ZZ$). 
\end{proposition} 
 
\proof 
Let us identify each cluster variable 
$x_i (t)$, for $t \in \TT_n$ and $i \in I_1 \subset I$ 
(resp., $i \in I_2 \subset I$), with $x_i (\pi_1 (t)) \otimes 1$ 
(resp., $1 \otimes x_i (\pi_2 (t))$. 
Under this identification, the exchange relations 
for the exchange pattern $\MM$ 
become identical to the exchange relations for $\MM_1$ and 
$\MM_2$. 
\endproof

\section{The Laurent phenomenon} 
\label{sec:laurent} 
 
In this section we prove the following important property of 
cluster algebras.

\begin{theorem} 
\label{th:laurent-binomial} 
In a cluster algebra, any cluster variable is expressed in terms of any 
given cluster as a Laurent polynomial with coefficients in the 
%integer 
group ring 
%of the coefficient group 
$\ZZ \PP$. 
\end{theorem} 
 
We conjecture that each of the coefficients in these Laurent polynomials 
is actually a \emph{nonnegative} integer linear combination of elements in~$\PP$. 
 
We will obtain Theorem~\ref{th:laurent-binomial} as a corollary of 
a more general result, which applies to %more general coefficient systems, 
more general underlying graphs and more general 
(not necessarily binomial) exchange polynomials. 
 
\begin{figure}[ht] 
\begin{center} 
\begin{picture}(280,60)(0,0) 
\thicklines 
  \put(0,30){\line(1,0){280}} 
  \multiput(20,30)(40,0){7}{\vector(1,0){7}} 
  \multiput(0,30)(40,0){8}{\line(1,-2){10}} 
  \multiput(0,30)(40,0){8}{\line(-1,-2){10}} 
  \put(0,30){\line(-2,1){20}} 
  \put(280,30){\line(2,1){20}} 
  \multiput(0,30)(40,0){8}{\circle*{4}} 
  \multiput(-10,10)(20,0){16}{\circle*{4}} 
  \put(-20,40){\circle*{4}} 
  \put(300,40){\circle*{4}} 
  \put(-30,47){$t_{\rm tail}$} 
  \put(-8,37){$t_{\rm base}$} 
  \put(290,47){$t_{\rm head}$} 
\end{picture} 
\end{center} 
\caption{The ``caterpillar'' tree $\TT_{n,m}$, for $n=4$, $m=8$} 
\label{fig:caterpillar} 
\end{figure}
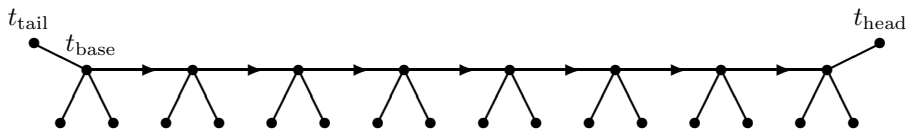

Since Theorem~\ref{th:laurent-binomial} is trivial for $n = 1$, 
we shall assume that $n \geq 2$. 
For every $m \geq 1$, let $\TT_{n,m}$ be a tree of the form shown in 
Figure~\ref{fig:caterpillar}. 
The tree $\TT_{n,m}$ has $m$ vertices of degree $n$ in its 
``spine" and $m (n-2)+2$ vertices of degree~1. 
%; all  vertices of degree $n$ lie on a path (``stem'') of length~$m$. 
We label every edge of the tree by an element of an $n$-element index set 
$I$, so that the $n$ edges incident to each vertex on the spine receive different 
labels. 
(The reader may wish to think of the tree $\TT_{n,m}$ as being part of 
%a larger tree, say 
the $n$-regular tree $\TT_n$ of the cluster-algebra setup.) 
 
We fix two vertices $t_{\rm head}$ and $t_{\rm tail}$ of $\TT_{n,m}$ that do not belong 
to the spine and are connected to its opposite ends. 
This gives rise to the orientation on the spine: away from $t_{\rm tail}$ and 
towards $t_{\rm head}$ (see Figure~\ref{fig:caterpillar}).

As before, let $\PP$ be an abelian group without torsion, written 
multiplicatively. 
Let $\ZZ_{\geq 0} \PP$ denote the additive semigroup generated by $\PP$ 
in the integer group ring $\ZZ \PP$. 
Assume that a nonzero polynomial $P$ in the variables $x_i \, i \in I$, 
with coefficients in $\ZZ_{\geq 0} \PP$, is associated with every edge 
$t\overunder{}{}t'$ of $\TT_{n,m}$. 
We call $P$ an \emph{exchange polynomial,} 
and write $t \overunder{}{P} t'$ to describe this situation. 
Suppose that the exchange polynomials 
associated with the edges of $\TT_{n,m}$ 
%is called a \emph{generalized exchange pattern} on~$\TT_{n,m}$ 
satisfy the following conditions: 
\begin{eqnarray} 
\label{eq:GEP1} 
&&\text{An exchange polynomial associated with an edge labeled by~$j$ 
  does not} 
\\ 
\nonumber 
&&\text{depend on $x_j$, and is not divisible by any $x_i$, $i\in I$.}\\ 
\label{eq:GEP2} 
&&\text{If $t_0 \overunder{i}{P} t_1 \overunderarrow{j}{Q} t_2 
\overunder{i}{R} t_3$, %(note that two edges are labelled $i$), 
then $R=C\cdot (P|_{x_j\leftarrow {Q_0}/{x_j}})$, where $Q_0=Q|_{x_i=0}\,$,}\\ 
\nonumber 
&&\text{and $C$ is a Laurent polynomial with coefficients in $\ZZ_{\geq 0} \PP$.} 
\end{eqnarray} 
\noindent (Note the orientation of the edge $t_1\to t_2$ in (\ref{eq:GEP2}).) 
 
\smallskip 
 
For every vertex $t$ on the spine, let $\mathcal{P}(t)$ 
denote the family of $n$ exchange polynomials 
associated with the edges emanating from~$t$. 
Also, let $\mathcal {C}$ denote the collection of all Laurent 
polynomials $C$ that appear in condition (\ref{eq:GEP2}), 
for all possible choices of $t_0, t_1, t_2, t_3$, and 
let $\AAA \subset \ZZ \PP$ denote the subring with unit generated by all coefficients 
of the Laurent polynomials from 
$\mathcal{P}(t_{\rm base}) \cup \mathcal {C}$, where $t_{\rm base}$ is the vertex on 
the spine connected with $t_{\rm tail}$. 
 
As before, 
%in Section~\ref{sec:intro}, 
we associate a 
\emph{cluster} ${\bf x}(t)=\{x_i(t)\ : \ i \in I\}$ to each 
vertex $t\in \TT_{n,m}\,$, 
and consider the field $\FFcal(t)$ of rational functions in these 
variables with coefficients in~$\ZZ \PP$. 
All these fields are identified with each other by the 
transition isomorphisms $R_{tt'}: \FFcal(t') \to \FFcal(t)$ 
defined as in (\ref{eq:transition1}). 
We then view the fields $\FFcal(t)$ as a single field $\FFcal$ that 
contains all the elements $x_i(t)$, for $t\in\TT_{n,m}$ and $i\in I$. 
These elements satisfy the exchange relations 
(\ref{eq:exchange1}) and the following version of (\ref{eq:exchange2}): 
$$x_j(t)\,x_j(t')= P({\bf x}(t))$$ 
for any edge $t \overunder{j}{P} t'$ in $\TT_{n,m}$.

\begin{theorem} 
\label{th:laurent-general} 
%For every generalized exchange pattern on~$\TT_{n,m}$, 
If conditions {\rm(\ref{eq:GEP1})--(\ref{eq:GEP2})} are satisfied, 
then each element of the cluster ${\bf x} (t_{\rm head})$ is a Laurent polynomial 
in the cluster ${\bf x}(t_{\rm tail})$, with coefficients in the ring~$\AAA$. 
%Properties {\rm(\ref{eq:GEP1})--(\ref{eq:GEP3})} imply that 
%whenever vertex $t'$ is reachable from $t$, 
%the transition map $R_{t't}$ is given by Laurent polynomials. 
\end{theorem} 
 
We note that Theorem~\ref{th:laurent-general} is indeed a generalization of 
Theorem~\ref{th:laurent-binomial}, for the following reasons: 
\begin{itemize} 
\item $\TT_{n,m}$ is naturally embedded into $\TT_n$; 
\item conditions  (\ref{eq:GEP1})--(\ref{eq:GEP2}) 
  are less restrictive than (\ref{eq:CA1})--(\ref{eq:CA4}); 
\item the claim being made in Theorem~\ref{th:laurent-general} 
about coefficients of the Laurent polynomials is stronger 
than that of Theorem~\ref{th:laurent-binomial}, since $\AAA \subset\ZZ\PP$. 
\end{itemize} 
 
\proof 
We start with some preparations. 
We shall write any Laurent polynomial $L$ in the variables 
${\bf x} = \{x_i : i \in I\}$ 
in the form 
$$L ({\bf x}) = \sum_{\alpha \in S} u_\alpha (L) x^\alpha \ ,$$ 
where all coefficients $u_\alpha (L)$ are nonzero, $S$ is a finite subset of the 
lattice $\ZZ^I$ (i.e., the lattice of rank $n$ with coordinates labeled by $I$), 
and $x^\alpha$ is the usual shorthand for 
$\prod_i x_i^{\alpha_i}$. 
The set $S$ is called the \emph{support} of $L$ and denoted by $S = S(L)$. 
 
Notice that once we fix the collection $\mathcal {C}$, condition (\ref{eq:GEP2}) can be 
used as a recursive rule for computing $\mathcal{P}(t')$ from $\mathcal{P}(t)$, for 
any edge $t \overunderarrow{j}{} t'$ on the spine. 
It follows that the whole pattern of exchange polynomials is determined by the 
families of polynomials $\mathcal{P}(t_{\rm base})$ and $\mathcal {C}$. 
Moreover, since these polynomials have coefficients in $\ZZ_{\geq 0} \PP$, 
and the expression for $R$ in (\ref{eq:GEP2}) does not involve subtraction, 
it follows that the support of any exchange polynomial is uniquely determined by 
the supports of the polynomials from $\mathcal{P}(t_{\rm base})$ and $\mathcal {C}$. 
Note that condition (\ref{eq:GEP1}) can be formulated as a set of 
restrictions on these supports. 
In particular, it requires that in the situation of 
(\ref{eq:GEP2}), the Laurent polynomial $C$ does not depend on 
$x_i$ and is a \emph{polynomial} in $x_j$; in other words, 
every $\alpha \in S(C)$ should have  $\alpha_i = 0$ and $\alpha_j \geq 0$. 
 
We now fix a family of supports $S(L)$, for all 
$L \in \mathcal{P}(t_{\rm base}) \cup \mathcal {C}$, and assume 
that this family complies with (\ref{eq:GEP1}). 
As is common in algebra, we shall view the 
coefficients $u_\alpha (L)$, for all 
$L \in \mathcal{P}(t_{\rm base}) \cup \mathcal {C}$ and 
$\alpha \in S(L)$, as \emph{indeterminates}. 
Then all the coefficients in all exchange polynomials 
become ``canonical" (i.e., independent of the choice of 
$\PP$) polynomials in these indeterminates, with 
positive integer coefficients. 
 
The above discussion shows that it suffices to prove our theorem in the 
following ``universal coefficients" setup: let $\PP$ be the free 
abelian group (written multiplicatively) with generators $u_\alpha (L)$, for 
all $L \in \mathcal{P}(t_{\rm base}) \cup \mathcal {C}$ and 
$\alpha \in S(L)$. 
Under this assumption, $\AAA$ is simply the integer 
polynomial ring in the indeterminates $u_\alpha (L)$. 
 
Recall that we can view all cluster variables $x_i (t)$ as elements of the 
field $\FFcal(t_{\rm tail})$ of rational functions in the cluster 
$x (t_{\rm tail})$ with coefficients in~$\ZZ \PP$. 
For $t \in \TT_{n,m}$, let $\LL (t)$ 
%\[ 
%\LL (t) = \AAA [x_1(t)^{\pm 1}, \dots, x_n(t)^{\pm 1}] 
%\] 
denote the ring of Laurent polynomials in the cluster ${\bf x}(t)$, with coefficients in~$\AAA$. 
We view each $\LL (t)$ as a subring of the ambient field $\FFcal(t_{\rm tail})$. 
 
In this terminology, our goal is to show that 
the cluster ${\bf x} (t_{\rm head})$ is contained in~$\LL (t_{\rm tail})$. 
We proceed by induction on $m$, the size of the spine. 
The claim is trivial for $m = 1$, so let us assume that 
$m \geq 2$, and furthermore assume that our statement is true for all 
``caterpillars" with smaller spine. 
 
Let us abbreviate $t_0 = t_{\rm tail}$ and $t_1 = t_{\rm base}$, 
and suppose that the path from $t_{\rm tail}$ to $t_{\rm head}$ starts 
with the following two edges: 
$t_0 \overunder{i}{P} t_1 \overunderarrow{j}{Q} t_2$. 
Let $t_3 \in \TT_{n,m}$ be the vertex such that 
$t_2 \overunder{i}{R} t_3$. 
 
The following lemma plays a crucial role in our proof. 
 
\begin{lemma} 
\label{lem:3-step-general} 
%Let 
%$~~t_0 \overunder{i}{P} t_1 \overunderarrow{j}{Q} t_2 \overunder{i}{R} t_3~$ 
%be a fragment of a generalized exchange pattern. 
%Then $x_i (t_3) \in \LL_0=\LL(t_0)$, and therefore all variables in the 
The clusters ${\bf x}(t_1)$, ${\bf x}(t_2)$, and ${\bf x}(t_3)$ are contained in 
$\LL (t_0)$. 
Furthermore, %we have 
$\gcd (x_i (t_3), x_i (t_1)) = \gcd (x_j (t_2), x_i (t_1)) = 1$ 
(as elements of $\LL (t_0)$). 
\end{lemma} 
 
Note that $\LL_0 = \LL(t_0)$ is a unique factorization domain, so 
any two elements $x, y \in \LL_0$ have a 
well-defined greatest common divisor $\gcd (x,y)$, 
which is an element of $\LL_0$ defined up to a multiple from the group 
$\LL_0^\times$ of units (that is, invertible elements) of $\LL_0$. 
In our ``universal'' situation, $\LL_0^\times$ consists of Laurent 
monomials in the cluster ${\bf x}(t_0)$ with coefficients $\pm 1$. 
 
\proof 
The only element from the clusters ${\bf x}(t_1)$, 
${\bf x}(t_2)$, and ${\bf x}(t_3)$ 
whose inclusion in $\LL_0$ is not immediately obvious is 
$x_i (t_3)$. 
To simplify the notation, let us denote 
$x=x_i(t_0)$, $y=x_j(t_0)=x_j(t_1)$, $z=x_i(t_1)=x_i(t_2)$, 
$u=x_j(t_2)=x_j(t_3)$, and $v=x_i(t_3)$, so that these variables 
appear in the clusters at $t_0,\dots,t_3$, as shown below: 
\[ 
\begin{array}{c}\scriptstyle x,y\\ \bullet\\ \scriptstyle t_0\end{array} 
 \hspace{-.1in}\wideoverunder{i}{P}\hspace{-.1in} 
\begin{array}{c}\scriptstyle y,z\\ \bullet\\ \scriptstyle t_1\end{array} 
 \hspace{-.1in}\wideoverunderarrow{j}{Q}\hspace{-.1in} 
\begin{array}{c}\scriptstyle z,u\\ \bullet\\ \scriptstyle t_2\end{array} 
 \hspace{-.1in}\wideoverunder{i}{R}\hspace{-.1in} 
\begin{array}{c}\scriptstyle u,v\\ \bullet\\ \scriptstyle t_3\end{array} 
\,. 
\] 
Note that the variables $x_k$, for $k\notin\{i,j\}$, do not change as 
we move among the four clusters under consideration. 
The lemma is then restated as saying that 
\begin{eqnarray} 
\label{eq:Laurent-3} 
&&\text{$v\in\LL_0$;}\\ 
\label{eq:gcd(u,z)} 
&&\text{$\gcd(z,u)=1$ (as elements of $\LL_0$);}\\ 
\label{eq:gcd(v,z)} 
&&\text{$\gcd(z,v)=1$ (as elements of $\LL_0$).} 
\end{eqnarray} 
Another notational convention will be based on the fact that each of 
the polynomials $P,Q,R$ has a distinguished variable on which it 
depends, namely $x_j$ for $P$ and $R$, and $x_i$ for~$Q$. 
(In view of (\ref{eq:GEP1}), $P$ and $R$ do not depend on $x_i$, while 
$Q$ does not depend on~$x_j$.) 
With this in mind, we will routinely write $P$, $Q$, and $R$ as 
polynomials in one (distinguished) variable. 
In the same spirit, the notation $Q'$, $R'$, etc., will refer to the 
partial derivative with respect to the distinguished variable. 
 
We will prove the statements (\ref{eq:Laurent-3}), 
(\ref{eq:gcd(u,z)}), and (\ref{eq:gcd(v,z)}) one by one, in this 
order. 
 
By (\ref{eq:GEP2}), the polynomial $R$ is given by 
\begin{equation} 
\label{eq:R(u)=C(u)P} 
R(u)=C(u)P\left(\textstyle\frac{Q(0)}{u}\right), 
\end{equation} 
where $C$ is an ``honest" polynomial in $u$ and a Laurent polynomial in 
the ``mute'' variables $x_k$, $k\notin\{i,j\}$. 
(Recall that $C$ does not depend on~$x_i$.) 
We then have: 
\begin{eqnarray*} 
&&z=\frac{P(y)}{x}\,;\\ 
&&u=\frac{Q(z)}{y}=\frac{Q\left(\frac{P(y)}{x}\right)}{y}\,;\\ 
&&v=\frac{R(u)}{z} 
=\frac{R\left(\frac{Q(z)}{y}\right)}{z} 
=\frac{R\left(\frac{Q(z)}{y}\right)-R\left(\frac{Q(0)}{y}\right)}{z} 
+\frac{R\left(\frac{Q(0)}{y}\right)}{z}\,. 
\end{eqnarray*} 
Since 
\[ 
\frac{R\left(\frac{Q(z)}{y}\right)-R\left(\frac{Q(0)}{y}\right)}{z} 
\in\LL_0 
\] 
and 
\[ 
\frac{R\left(\frac{Q(0)}{y}\right)}{z} 
=\frac{C \left(\frac{Q(0)}{y}\right)P(y)}{z} 
=C \left(\textstyle\frac{Q(0)}{y}\right)x\in\LL_0 \,, 
\] 
(\ref{eq:Laurent-3}) follows. 
 
We next prove (\ref{eq:gcd(u,z)}). 
We have 
\[ 
u=\frac{Q(z)}{y}\equiv \frac{Q(0)}{y}\bmod z\,. 
\] 
%If $Q(z)$ depends on $z$, then by (\ref{eq:CA2e}), 
%$\frac{Q(0)}{y}$ is an invertible element in $\LL_0$, 
%and there is nothing to prove. 
%If, on the other hand, $Q(z)$ does not depend on $z$, then 
Since $x$ and $y$ are invertible in $\LL_0$, we conclude that 
$\gcd(z,u)=\gcd(P(y),Q(0))$. 
Now the trouble that we took in passing to universal 
coefficients finally pays off: since $P(y)$ and $Q(0)$ 
are nonzero polynomials in the cluster ${\bf x}(t_0)$ 
whose coefficients are distinct generators of 
the polynomial ring $\AAA$, it follows that 
$\gcd(P(y),Q(0)) = 1$, proving (\ref{eq:gcd(u,z)}). 
%(Note that (\ref{eq:GEP1}) implies $Q(0)\neq 0$.) 
 
%Let us prove (\ref{eq:gcd(v,z)}). 
%If $Q(z)$ does not depend on $z$, then by (\ref{eq:CA3e}), 
%$R$ and $P$ do not depend on their distinguished variables either. 
%Moreover, (\ref{eq:CA4e}) implies that in that case, $P$ and $Q$ 
%differ by a constant factor. 
%Hence $v=\frac{R(u)}{z}\multsim\frac{P(y)}{z}=x$, 
%and (\ref{eq:gcd(v,z)}) is trivial. 
 
It remains to prove (\ref{eq:gcd(v,z)}). 
%under the assumption that $Q(z)$ does depend on~$z$. 
Let 
\[ 
f(z)=R\left(\textstyle\frac{Q(z)}{y}\right). 
\] 
Then 
\[ 
v=\frac{f(z)-f(0)}{z}+C \left(\textstyle\frac{Q(0)}{y}\right)x\,. 
\] 
Our goal is to show that $\gcd(z,v)=1$; to this end, we are going to 
compute $v\bmod z$ as ``explicitly" as possible. 
We have, $\bmod z$, 
\[ 
\frac{f(z)-f(0)}{z} 
\equiv f'(0) 
=R'\left(\textstyle\frac{Q(0)}{y}\right) 
\cdot\textstyle\frac{Q'(0)}{y}\,. 
\] 
%%% 
%%%  From (\ref{eq:R(u)=A(u)P}) we get 
%%% \[ 
%%% R'(u) 
%%% =A'(u) \, P\!\left(\textstyle\frac{Q(0)}{u}\right) 
%%% -A(u)\, \frac{Q(0)}{u^2}\, P'\!\left(\textstyle\frac{Q(0)}{u}\right) , 
%%% \] 
%%% implying 
%%% \begin{eqnarray*} 
%%% R'\!\left(\textstyle\frac{Q(0)}{y}\right) 
%%% &=&A'\!\left(\textstyle\frac{Q(0)}{y}\right) \, P(y) 
%%% -A\!\left(\textstyle\frac{Q(0)}{y}\right)\, \frac{y^2}{Q(0)} 
%%% \, P'(y) \\ 
%%% &\equiv& -A\!\left(\textstyle\frac{Q(0)}{y}\right)\, \frac{y^2}{Q(0)} 
%%% \, P'(y)~ \bmod z\,. 
%%% \end{eqnarray*} 
%%% %(since $P(y)\equiv 0\bmod z$). 
%%% We conclude that, $\bmod z$, 
%%% \[ 
%%% v\equiv -A\!\left(\textstyle\frac{Q(0)}{y}\right)\, \frac{y Q'(0)}{Q(0)} 
%%% \, P'(y) 
%%% +A\left(\textstyle\frac{Q(0)}{y}\right)x\,. 
%%% \] 
Hence 
\[ 
v\equiv R'\left(\textstyle\frac{Q(0)}{y}\right) 
\cdot\textstyle\frac{Q'(0)}{y} 
+C \left(\textstyle\frac{Q(0)}{y}\right)x 
\bmod z\,. 
\] 
Note that the right-hand side is a \emph{linear} polynomial in $x$, 
whose coefficients are Laurent polynomials in the rest of the 
variables of the cluster ${\bf x}(t_0)$. 
Thus our claim will follow if we show that 
$\gcd\left(C \left(\textstyle\frac{Q(0)}{y}\right),P(y)\right)=1$. 
This, again, is a consequence of our ``universal coefficients" 
setup since the coefficients of $C$, $P$ and $Q$ are distinct 
generators of the polynomial ring~$\AAA$. 
\endproof

We can now complete the proof of Theorem~\ref{th:laurent-general}. 
We need to show that %for any  vertices $t_0,t\in T$, 
any variable $x = x_k(t_{\rm head})$ %, for $t\in T$ and $k\in[n]$, 
belongs to $\LL (t_0)$. 
Since both $t_1$ and $t_3$ are closer to $t_{\rm head}$ than 
$t_0$, we can use the inductive assumption to conclude 
that $x$ belongs to both $\LL(t_1)$ and $\LL(t_3)$. 
Since $x \in \LL(t_1)$, it follows from (\ref{eq:exchange1})  that 
$x$ can be written as $x = f/ x_i(t_1)^a$ for some $f \in \LL(t_0)$ 
and $a\in\ZZ_{\geq 0}\,$. 
On the other hand, since $x \in \LL(t_3)$, it follows from 
(\ref{eq:exchange1}) and from the inclusion $x_i (t_3) \in \LL(t_0)$ 
guaranteed by Lemma~\ref{lem:3-step-general} 
that $x$ has the form  $x = g/ x_j(t_2)^b x_i(t_3)^c$ for some $g \in \LL(t_0)$ 
and some $b, c\in\ZZ_{\geq 0}\,$. 
The inclusion $x \in \LL(t_0)$ now follows from the fact 
that, by the last statement in Lemma~\ref{lem:3-step-general}, 
the denominators in the two obtained expressions for $x$ are 
coprime in~$\LL (t_0)$. 
\endproof 
 
Several examples that can be viewed as applications of Theorem~\ref{th:laurent-general} 
are given in~\cite{fz-Laurent}.

\section{Exchange relations: the exponents} 
\label{sec:exponents} 
 
Let $\MM=(M_j (t)) \ : \ t\in\TT_n, j \in I)$ be an exchange pattern 
(see Definition~\ref{def:exchange pattern}). 
In this section we will ignore the 
coefficients in the monomials $M_j (t)$ and take a closer look at 
the dynamics of their exponents. 
(An alternative point of view that the reader may find helpful is to assume 
that all exchange patterns considered in this section will have
all their coefficients $p_j (t)$ equal to~$1$.) 
For every edge $t\overunder{j}{}t'$ in $\TT_n$, let us write the 
ratio $M_j(t)/M_j(t')$ of the corresponding monomials as 
\begin{equation} 
\label{eq:B(t)} 
\frac{M_j(t)}{M_j(t')}=\frac{p_j(t)}{p_j(t')} 
\prod_{i} x_i^{b_{ij}(t)} \ , 
\end{equation} 
where $b_{ij} (t) \in \ZZ$ (cf.~(\ref{eq:M(t)})); we note that ratios of this kind 
have already appeared in (\ref{eq:CA4}).  
Let us denote by $B(t)=(b_{ij}(t))$ the $n\times n$ integer matrix 
whose entries are the exponents in (\ref{eq:B(t)}). 
In view of (\ref{eq:CA2}), the exponents in $M_j(t)$ and 
$M_j(t')$ are recovered from $B(t)$: 
\begin{equation} 
\label{eq:M_k(t)} 
\begin{array}{rcl} 
M_j(t) & = & p_j(t) \displaystyle\prod_{i: \ b_{ij}(t)>0} 
x_i^{b_{ij}(t)}\, , 
\\[.2in] 
M_j(t')&=& p_k(t') \displaystyle\prod_{i: \ b_{ij}(t) <0} x_i^{-b_{ij}(t)} 
\ . 
%\qquad\text{(here $~t\overunder{k}{}t'$)}, 
\end{array} 
\end{equation} 
Thus, the family of matrices $(B(t))_{t \in \TT_n}$ encodes all the exponents in 
all monomials of an exchange pattern. 
 
We shall describe the conditions on the family of matrices $(B(t))$ imposed by the axioms 
of an exchange pattern. 
To do this, we need some preparation. 
 
\begin{definition} 
\label{def:sign-ss} 
{\rm A square integer matrix $B = (b_{ij})$ is called 
\emph{sign-skew-symmetric} if, for any $i$ and $j$, 
either $b_{ij} = b_{ji} = 0$, or else $b_{ij}$ and $b_{ji}$ 
are of opposite sign; in particular, $b_{ii} = 0$ for all $i$.} 
\end{definition} 
 
\begin{definition} 
\label{def:matrix mutation} 
{\rm Let $B = (b_{ij})$ and $B' = (b'_{ij})$ be square integer matrices of the same size. 
We say that $B'$ is obtained from $B$ by the \emph{matrix mutation} in direction~$k$ 
and write $B' = \mu_k (B)$ if 
\begin{equation} 
\label{eq:mutation} 
b'_{ij} = 
\begin{cases} 
-b_{ij} & \text{if $i=k$ or $j=k$;} \\[.15in] 
b_{ij} + \displaystyle\frac{|b_{ik}| b_{kj} + 
b_{ik} |b_{kj}|}{2} & \text{otherwise.} 
\end{cases} 
\end{equation} 
} 
\end{definition} 
 
An immediate check shows that $\mu_k$ is \emph{involutive}, i.e., 
its square is the identity transformation. 
 
\begin{proposition} 
\label{prop:B(t)} 
A family of $n\times n$ integer matrices %$B(t) = (b_{ij}(t))_{t\in\TT_n}$ 
$(B(t))_{t\in\TT_n}$ corresponds to an exchange pattern 
if and only if the following conditions hold:\\ 
(1) $B(t)$ is sign-skew-symmetric for any  $t\in\TT_n$.\\
(2) If $~t \!\overunder{k}{}t'$, then $B(t') = \mu_k (B(t))$. 
\end{proposition}

\proof 
We start with the ``only if" part, i.e., we assume that the matrices 
$B(t)$ are determined by an exchange pattern via (\ref{eq:B(t)}) and check the conditions  
(1)--(2). 
The condition $b_{jj} (t) = 0$ follows from (\ref{eq:CA1}). 
The remaining part of (1) (dealing with $i \neq j$), 
follows at once from (\ref{eq:CA3}). 
Turning to part (2), the equality $b'_{ik}  = - b_{ik}$ is 
immediate from the definition (\ref{eq:B(t)}). 
Now suppose that $j \neq k$. 
In this case, we apply the axiom (\ref{eq:CA4}) to the edge 
$t \!\overunder{k}{}t'$ taken together with the two adjacent edges emanating 
from $t$ and $t'$ and labeled by $j$. 
Taking (\ref{eq:M_k(t)}) into account, we obtain: 
\begin{equation} 
\label{eq:CA4 thru B} 
\prod_{i} x_i^{b'_{ij}}= \prod_{i} x_i^{b_{ij}} \Bigl|_{x_k\leftarrow 
{M}/{x_k}} \ , 
\end{equation} 
where $M = \displaystyle \prod_{i: b_{ik} b_{jk} < 0} x_i^{|b_{ik}|}$. 
Comparing the exponents of $x_k$ on both sides of 
(\ref{eq:CA4 thru B}) yields $b'_{kj} = - b_{kj}$. 
Finally, if $i \neq k$ then, comparing the exponents of $x_i$ on both sides of 
(\ref{eq:CA4 thru B}) gives 
$$b'_{ij} = 
\begin{cases} 
b_{ij} & \text{if $b_{ik} b_{jk} \geq 0$ ;} \\[.15in] 
b_{ij} + |b_{ik}| b_{kj}  & \text{otherwise.} 
\end{cases}$$ 
To complete the proof of (2), it remains to notice that, 
in view of the already proven part (1), the condition 
$b_{ik} b_{jk} \geq 0$ is equivalent to $b_{ik} b_{kj} \leq 0$, which 
makes the last formula equivalent to (\ref{eq:mutation}). 
 
To prove the ``if" part, it suffices to show that if the matrices 
$B(t)$ satisfy (1)--(2) then the monomials 
$M_j (t)$ given by the first equality in (\ref{eq:M_k(t)}) 
(with $p_j (t) = 1$) satisfy the axioms of an exchange pattern. 
This is done by a direct check. \endproof

Since all matrix mutations are involutive, 
any choice of an initial vertex $t_0 \in \TT_n$ 
and an arbitrary $n \times n$ integer matrix $B$ gives rise to a 
unique family of integer matrices $B(t)$ satisfying condition (2) in 
Proposition~\ref{prop:B(t)} and such that $B(t_0) = B$. 
Thus, the exponents in all monomials $M_j (t)$ are uniquely 
determined by a single matrix $B = B(t_0)$. 
By Proposition~\ref{prop:B(t)}, in order to determine an exchange 
pattern, $B$ must be such that all matrices obtained 
from it by a sequence of matrix mutations are sign-skew-symmetric. 
Verifying that a given matrix $B$ has this property 
seems to be quite non-trivial in general. 
Fortunately, there is another restriction on $B$ that is much easier to check, 
which implies the desired property, and still leaves us with 
a large class of matrices sufficient for most applications. 
 
\begin{definition} 
{\rm 
A square integer matrix $B = (b_{ij})$ is called 
\emph{skew-symmetrizable} if there exists a 
diagonal \emph{skew-symmetrizing} matrix $D$ with 
positive integer diagonal entries $d_i$ such that $DB$ is 
skew-symmetric, i.e., $d_i b_{ij} = - d_j b_{ji}$ 
for all $i$ and $j$. 
%$(DB)^{\rm transpose}=-DB$. 
} 
\end{definition} 
 
%Note that every skew-symmetrizable matrix $B$ satisfies 
%(\ref{eq:mock-skew-symm}). 
 
%The following simple result explains the significance of 
%skew-symmetrizability. 
 
\begin{proposition} 
\label{prop:skew-symm-from-t} 
For every choice of a vertex $t_0 \in \TT_n$ and a 
skew-symmetrizable matrix~$B$, there exists a unique family 
of matrices $(B(t))_{t \in \TT_n}$ associated with an exchange 
pattern on $\TT_n$ and such that $B(t_0) = B$. 
Furthermore, all the matrices $B(t)$ are 
skew-symmetrizable, sharing the same skew-symmetrizing matrix. 
\end{proposition} 
 
\proof 
The proof follows at once from the following two observations: \\ 
1. Every skew-symmetrizable matrix $B$ is sign-skew-symmetric. \\
2. If $B$ is skew-symmetrizable, and $B' = \mu_k (B)$ 
then $B'$ is also skew-symmetrizable, with the same skew-symmetrizing matrix. 
\endproof 
 
We call an exchange pattern---and the corresponding cluster 
algebra---\emph{skew-symmetrizable} if all the matrices $B(t)$ given 
by (\ref{eq:B(t)}) (equivalently, one of them) are skew-symmetrizable. 
In particular, all cluster algebras of rank $n \leq 2$ are 
skew-symmetrizable: for $n = 1$ we have $B (t) \equiv (0)$, while 
for $n = 2$, the calculations in Example~\ref{ex:rank 2} show that 
one can take 
\begin{equation} 
\label{eq:B(t)-rank2} 
B(t_m)= (-1)^{m} \bmat{0}{b}{-c}{0} 
\end{equation} 
for all $m \in \ZZ$, in the notation of (\ref{eq:tree-2})--(\ref{eq:exchange rank 2}). 
 
\begin{remark} 
\label{rem:Kac-Moody connection} 
{\rm 
Skew-symmetrizable matrices are closely related to symmetrizable 
(generalized) Cartan matrices appearing in the theory of Kac-Moody algebras. 
More generally, to every sign-skew-symmetric matrix $B=(b_{ij})$ 
%satisfying condition (1) in Proposition~\ref{prop:B(t)} 
we can associate a generalized Cartan matrix $A = A(B)=(a_{ij})$ of the same size by setting 
\begin{equation} 
\label{eq:cartan matrix} 
a_{ij} = 
\begin{cases} 
2 & \text{if $i=j$;} \\[.1in] 
- |b_{ij}| & \text{if $i\neq j$.} 
\end{cases} 
\end{equation} 
There seem to be deep connections between the cluster algebra 
corresponding to $B$ and the Kac-Moody algebra associated with $A(B)$. 
We exhibit such a connection for the rank $2$ case in Section~\ref{sec:rank 2} below. 
This is however just the tip of an iceberg: a much more detailed 
analysis will be presented in the sequel to this paper.} 
\end{remark} 
 
In order to show that non-skew-symmetrizable exchange patterns do 
exist, we conclude this section by exhibiting a 3-parameter 
family of such patterns of rank~$3$. 
 
\begin{proposition} 
\label{non-ss rank 3} 
Let $\alpha$, $\beta$, and $\gamma$ be three positive integers such 
that $\alpha \beta \gamma \geq 3$. 
There exists a unique family 
of matrices $(B(t))_{t \in \TT_3}$ associated with a 
non-skew-symmetrizable exchange 
pattern on $\TT_3$ and such that the matrix $B(t_0)$ at a given 
vertex $t_0 \in \TT_3$ is equal to 
\begin{equation} 
\label{eq:non-ss} 
B(\alpha, \beta, \gamma) = 
\left[\begin{array}{ccc} 
0    &  2 \alpha         & - 2 \alpha \beta \\ 
- \beta \gamma  & 0 & 2 \beta \\ 
\gamma  & - \alpha \gamma & 0 \\ 
\end{array}\right] . 
\end{equation} 
\end{proposition}

\proof 
First of all, the matrix $B(\alpha, \beta, \gamma)$ is 
sign-skew-symmetric but not skew-symmetrizable. 
Indeed, any skew-symmetrizable matrix 
$B = (b_{ij})$ satisfies the equation 
$b_{12} b_{23} b_{31} = - b_{21} b_{32} b_{13}$. 
However, this equation for $B(\alpha, \beta, \gamma)$ holds only 
when $\alpha \beta \gamma$ is equal to $0$ or $2$. 
 
For the purpose of this proof only, we refer to a $3 \times 3$ 
matrix $B$ as \emph{cyclical} if its entries follow one of the two 
sign patterns 
$$\left[\begin{array}{ccc} 
0    &  +  & - \\ 
-  & 0 & +\\ 
+  & -  & 0 \\ 
\end{array}\right] 
\ ,  
\left[\begin{array}{ccc} 
0    &  -  & +\\ 
+  & 0 & -\\ 
-  & +  & 0 \\ 
\end{array}\right] 
\ .$$ 
In particular, $B(\alpha, \beta, \gamma)$ is cyclical; 
to prove the proposition, it suffices to show that any matrix 
obtained from it by a sequence of matrix mutations is also 
cyclical. 
 
The set of cyclical matrices is not stable under matrix mutations. 
Let us define some subsets of cyclical matrices that behave 
nicely with respect to matrix mutations. 
For a $3 \times 3$ matrix $B$, we denote 
$$c_1 = |b_{23} b_{32}|, \,\, c_2 = |b_{13} b_{31}|, \,\, 
c_3 = |b_{12} b_{21}|, \,\, r = |b_{12} b_{23} b_{31}| \ .$$ 
For $i \in \{1,2,3\}$, we say that $B$ is \emph{$i$-biased} if we have 
$$r > c_i \geq r/2 \geq c_j \geq 6 $$ 
for any $j \in \{1,2,3\} \setminus \{i\}$. 
For $B(\alpha, \beta, \gamma)$, we have 
$$r/2 = c_1 = c_2 = c_3 = 2 \alpha \beta \gamma \geq 6 \ ,$$ 
so it is $i$-biased for every $i \in \{1,2,3\}$. 
 
Our proposition becomes an immediate consequence of the following lemma. 
 
\begin{lemma} 
\label{lem:bias} 
Suppose $B$ is cyclical and $i$-biased, and let $j \neq i$. 
Then $\mu_j (B)$ is cyclical and $j$-biased. 
\end{lemma} 
 
\proof 
Without loss of generality we can assume that $i = 1$ and $j = 2$. 
Denote $B' = \mu_2 (B)$, and let us write 
$r' =  |b'_{12} b'_{23} b'_{31}|$, $c'_1 = |b'_{23} b'_{32}|$, etc. 
By (\ref{eq:mutation}), we have $b'_{12} = - b_{12}, \, b'_{21} = - b_{21}, \, 
b'_{23} = - b_{23}, \, b'_{32} = - b_{32} \ ;$ 
therefore, $c'_1 = c_1$ and $c'_3 = c_3$. 
We also have 
\begin{equation} 
\label{eq:b13new} 
b'_{13} = b_{13} + \frac{|b_{12}| b_{23} + 
b_{12} |b_{23}|}{2}, \, \, 
b'_{31} = b_{31} + \frac{|b_{32}| b_{21} + 
b_{32} |b_{21}|}{2} \ . 
\end{equation} 
Since $B$ is cyclical, the two summands on the right-hand side of each of the equalities in 
(\ref{eq:b13new}) have opposite signs. 
Note that 
$$|b_{12} b_{23}|  = \frac{r}{|b_{31}|} \geq \frac{2 c_2}{|b_{31}|} = 
2 |b_{13}| \ ,$$ 
and 
$$|b_{32} b_{21}| = \frac{c_1 |b_{21}|}{|b_{23}|} \geq 
\frac{r |b_{21}|}{2|b_{23}|} = \frac{c_3 |b_{31}|}{2} \geq 3 |b_{31}| \ .$$ 
It follows that $b'_{13}$ (resp., $b'_{31}$) has the opposite sign to $b_{13}$ 
(resp., $b_{31}$). 
Thus, $B'$ is cyclical, and it only remains to show 
that $B'$ is $2$-biased. 
To this effect, we note that $|b'_{31}| = |b_{32} b_{21}| - |b_{31}|$, and so 
$$r' = |b'_{12} b'_{23} b'_{31}| = 
|b_{12}| |b_{23}| (|b_{32} b_{21}| - |b_{31}|) = c_1 c_3 - r 
\geq (r/2) \cdot 6 - r = 2r \ ;$$ 
therefore, both $c'_1 = c_1$ and $c'_3 = c_3$ do not exceed $r'/2$. 
As for $c'_2$, we have $|b'_{13}| = |b_{12} b_{23}| - |b_{13}|$, and so 
$$c'_2 = |b'_{13} b'_{31}| = ( |b_{12} b_{23}| - |b_{13}|) 
(|b_{32} b_{21}| - |b_{31}|) = r' (1 - \frac{c_2}{r}) \ ;$$ 
since $c_2/r \leq 1/2$, we conclude that 
$r' > r' (1 - \frac{c_2}{r}) = c'_2 \geq r'/2$. 
This completes the proof that $B'$ is $2$-biased. 
Lemma~\ref{lem:bias} and Proposition~\ref{non-ss rank 3} are proved. 
\endproof

\section{Exchange relations: the coefficients} 
\label{sec:coefficients} 
 
In this section we fix a family of matrices $B(t)$ satisfying 
the conditions in Proposition~\ref{prop:B(t)}, 
and discuss possible choices of coefficients $p_j (t)$ 
that can appear in the corresponding exchange pattern. 
We start with the following simple characterization. 
 
\begin{proposition} 
\label{pr:ppp=ppp} 
Assume that matrices $(B(t))_{t \in \TT_n}$ satisfy conditions {\rm (1)--(2)} 
in Proposition~\ref{prop:B(t)}. 
A family of elements $p_j (t)$ of a coefficient group 
$\PP$ gives rise, via {\rm (\ref{eq:M_k(t)})}, to an exchange pattern 
if and only if they satisfy the following relations, whenever 
$t_1 \overunder{i}{} t_2 \overunder{j}{} t_3 \overunder{i}{} t_4$: 
\begin{equation} 
\label{eq:relation-coefficients} 
p_i (t_1) p_i (t_3)  p_j (t_3)^{\max (b_{ji} (t_3),0)} = p_i (t_2) 
p_i (t_4) p_j (t_2)^{\max (b_{ji} (t_2),0)} \ . 
\end{equation} 
\end{proposition} 
 
Note that $b_{ji} (t_2) = - b_{ji} (t_3)$ by (\ref{eq:mutation}), so at most one of the 
$p_j (t_2)$ and $p_j (t_3)$ actually enters (\ref{eq:relation-coefficients}). 
 
\proof 
The only axiom of an exchange pattern that involves the 
coefficients is (\ref{eq:CA4}), and the relation 
(\ref{eq:relation-coefficients}) is precisely what it prescribes. 
\endproof 
 
First of all, let us mention the trivial solution of (\ref{eq:relation-coefficients}) 
when all the coefficients $p_i (t)$ are equal to $1$. 
 
Moving in the opposite direction, we 
introduce the \emph{universal coefficient group}  
(with respect to a fixed family $(B(t))$) 
as the abelian group $\PP$ generated by the elements $p_i (t)$, for all 
$i \in I$ and $t \in \TT_n$, which has (\ref{eq:relation-coefficients}) 
as the system of defining relations. 
The torsion-freeness of this group is guaranteed 
by the following proposition whose straightforward proof is omitted.

\begin{proposition} 
\label{pr:free generators universal} 
The universal coefficient group $\PP$ is a free abelian group. 
More specifically, let $t_0 \in \TT_n$, and let $S$ be a collection of 
pairs $(i,t)$ %with the following property: 
that contains both $(i,t_0)$ and $(i,t)$ 
for any edge $t_0 \overunder{i}{} t$, 
and precisely one of the pairs $(i,t)$ and $(i,t')$ 
for each edge $t \overunder{i}{} t'$ with $t$ and $t'$ 
different from $t_0$. 
Then $\{p_i(t): (i,t) \in S\}$ is a set of free generators for $\PP$. 
\end{proposition} 
 
We see that, in contrast to (\ref{eq:mutation}), relations 
(\ref{eq:relation-coefficients}) leave infinitely many degrees of 
freedom in determining the coefficients $p_i (t)$ 
(cf.~Remark~\ref{rem:degree of freedom}). 
 
The rest of this section is devoted to  
important classes of exchange patterns within which all the 
coefficients are completely determined by specifying $2n$ of 
them corresponding to the edges emanating from a given vertex. 
 
Suppose that, in addition to the multiplicative group structure, 
the coefficient group $\PP$ is supplied with a binary operation 
$\oplus$ that we call \emph{auxiliary addition}. 
Furthermore, suppose that this operation is commutative, associative and 
distributive with respect to multiplication; thus 
$(\PP,\oplus, \cdot)$ is a \emph{semifield}. 
(By the way, under these assumptions $\PP$ is automatically 
torsion-free as a multiplicative group: indeed, if $p^m = 1$ for some 
$p \in \PP$ and $m \geq 2$, then 
$$p = \frac{p^m \oplus p^{m-1} \oplus \cdots \oplus p} 
{p^{m-1} \oplus p^{m-2} \oplus \cdots \oplus 1} = 
\frac{p^{m-1} \oplus p^{m-2} \oplus \cdots \oplus 1} 
{p^{m-1} \oplus p^{m-2} \oplus \cdots \oplus 1} = 1 \ . )$$ 
 
\begin{definition} 
\label{def:normalized} 
{\rm 
An exchange pattern and the corresponding cluster algebra are 
called \emph{normalized} if $\PP$ is a semifield, and 
$p_j (t) \oplus p_j (t') = 1$ 
for any edge $t \overunder{j}{} t'$. 
} 
\end{definition}

\begin{proposition} 
\label{pr:normalized} 
Fix a vertex $t_0 \in \TT_n$, and $2n$ elements $q_j$ and 
$r_j$ ($j \in I$)  of a semifield $\PP$ such that $q_j \oplus r_j = 1$ for all $j$. 
Then every family of matrices $B(t)$ satisfying the conditions 
in Proposition~\ref{prop:B(t)} 
gives rise to a unique normalized exchange pattern such that, 
for every edge $t_0 \overunder{j}{} t$, we have $p_j (t_0) = q_j$ and $ p_j (t) = r_j$. 
Thus a normalized exchange pattern is completely determined by the $2n$ monomials
$M_j (t_0)$ and $M_j (t)$, for all edges $t_0 \overunder{j}{} t$. 
\end{proposition} 
 
\proof 
In a normalized exchange pattern, the coefficients 
$p_j (t)$ and $p_j (t')$ corresponding to an edge 
$t \overunder{j}{} t'$ are determined by their ratio 
\begin{equation}
u_j (t) = \frac{p_j (t)}{p_j (t')}
\end{equation}
via 
\begin{equation} 
\label{eq:p through u} 
p_j (t) = \frac{u_j (t)}{1 \oplus u_j (t)}, \quad 
p_j (t') = \frac{1}{1 \oplus u_j (t)} \ . 
\end{equation} 
Clearly, we have 
\begin{equation} 
\label{eq:uu'} 
u_j (t) u_j (t') = 1 
\end{equation} 
for any edge $t \overunder{j}{} t'$. 
We can also rewrite the relation (\ref{eq:relation-coefficients}) as follows: 
\begin{equation} 
\label{eq:3 relation u} 
u_i (t') =  u_i (t) u_j (t)^{\max(b_{ji} (t),0)} 
(1 \oplus u_j (t))^{- b_{ji} (t)} 
\end{equation} 
for any edge $t \overunder{j}{} t'$ and any $i \neq j$. 
This form of the relations for the normalized coefficients makes 
our proposition obvious. 
\endproof 
 
\begin{remark} 
\label{rem:extra relations} 
{\rm It is natural to ask whether the normalization condition imposes additional 
multiplicative relations among the coefficients $p_i (t)$ 
that are not consequences of (\ref{eq:relation-coefficients}). 
In other words: can a normalized 
system of coefficients generate the universal coefficient group? 
In the next section, we present a complete answer to this 
question in the rank $2$ case (see Remark~\ref{rem:normalized coefficient group} 
below).} 
\end{remark}

One example of a semifield is the 
multiplicative group $\RR_{> 0}$ of positive real numbers, 
$\oplus$ being ordinary addition. 
However the following example is more important for our 
purposes. 
 
\begin{example} 
\label{ex:coefficients-tropical} 
{\rm Let $\PP$ be a free abelian group, written multiplicatively, 
with a finite set of generators $p_i \, (i \in I')$, and with 
auxiliary addition $\oplus$ defined by 
 
\begin{equation} 
\label{eq:tropical addition} 
\prod_i p_i^{a_i} \oplus \prod_i p_i^{b_i}  = 
\prod_i p_i^{\min (a_i, b_i)}  \ . 
\end{equation} 
Then $\PP$ is a semifield; specifically, it is a product of $|I'|$ copies 
of the \emph{tropical} semifield (see, e.g., \cite{bfz96}). 
We denote this semifield by $\Trop (p_i: i \in I')$. 
Note~that if all exponents $a_i$ and $b_i$ in (\ref{eq:tropical addition}) 
are nonnegative, then the monomial on the right-hand side 
is the gcd of the two monomials on the left.
} 
\end{example}

\begin{definition} 
\label{def:geometric type} 
{\rm We say that a normalized exchange pattern is of \emph{geometric type} 
if $\PP = \Trop (p_i: i \in I')$, and each coefficient $p_j (t)$ is a monomial in the 
generators $p_i$ with all exponents \emph{nonnegative}. 
} 
\end{definition} 
 
For an exchange pattern of geometric type, the 
normalization condition 
%$p_i (t) \oplus p_i (t') = 1$ 
simply means that, for every edge $t \overunder{j}{} t'$, 
the two monomials $p_j (t)$ and $p_j (t')$ in the generators 
$p_i$ are coprime, i.e., have no variable in common. 
 
In all our examples of cluster algebras of geometric origin, 
including those discussed in the introduction, the exchange 
patterns turn out to be of geometric type. 
These patterns have the following useful equivalent description. 
 
\begin{proposition} 
\label{pr:geometric type} 
Let $\PP = \Trop (p_i: i \in I')$. 
A family of coefficients $p_j (t) \in \PP$ gives rise to an 
exchange pattern of geometric type if and only if 
they are given by 
\begin{equation} 
\label{eq:geometric coefficients} 
p_j (t) = \prod_{i \in I'} p_i^{\max (c_{ij} (t), 0)} 
\end{equation} 
for some family of integers $(c_{ij}(t): t \in \TT_n, i \in I', j \in I)$) 
satisfying the following property: 
for every edge $t \overunder{k}{} t'$ in $\TT_n$, 
the matrices $C(t) = (c_{ij} (t)) = (c_{ij})$ and $C(t') = (c_{ij}(t')) = (c'_{ij})$ 
are related by 
\begin{equation} 
\label{eq:C-mutation} 
c'_{ij} = 
\begin{cases} 
-c_{ij} & \text{if $j=k$;} \\[.15in] 
c_{ij} + \displaystyle\frac{|c_{ik}| b_{kj}(t) + 
c_{ik} |b_{kj}(t)|}{2} & \text{otherwise.} 
\end{cases} 
\end{equation} 
\end{proposition} 
 
\proof 
As in the proof of Proposition~\ref{pr:normalized}, 
for every edge $t \overunder{j}{} t'$ we consider the ratio 
$u_j (t) = p_j (t)/p_j (t')$. 
We introduce the matrices $C(t)$ by setting 
$$u_j (t) = \prod_i p_i^{c_{ij} (t)} \ .$$ 
The expression (\ref{eq:geometric coefficients}) then becomes 
a specialization of the first equality in (\ref{eq:p through u}), 
with auxiliary addition given by (\ref{eq:tropical addition}). 
To derive (\ref{eq:C-mutation}) from (\ref{eq:3 relation u}), 
first replace $i$ by $j$ and $j$ by $k$, respectively, then specialize,
then pick up the exponent of~$p_i\,$.  
\endproof 
 
Comparing (\ref{eq:C-mutation}) with (\ref{eq:mutation}), 
we see that it is natural to combine a pair of matrices 
$(B(t), C(t))$ into one rectangular integer matrix 
$\tilde B(t) = (b_{ij})_{i \in I \coprod I', j \in I}$ 
by setting $b_{ij} = c_{ij}$ for  $i \in I'$ and $j \in I$. 
Then the matrices $\tilde B (t)$ for $t \in \TT_n$ are related to 
each other by the same matrix mutation rule (\ref{eq:mutation}), 
now applied to any $i \in I \coprod I'$ and $j \in I$. 
We refer to $B(t)$ as the \emph{principal part} of $\tilde B (t)$. 
Combining Propositions~\ref{pr:geometric type} and 
\ref{prop:skew-symm-from-t}, we obtain the following corollary. 
 
\begin{corollary} 
\label{cor:geometric type local} 
Suppose $\tilde B$ is an integer matrix whose 
principal part is skew-symmetrizable. 
There is a unique exchange pattern $\MM = \MM(\tilde B)$ of geometric type such that 
$\tilde B (t_0) = \tilde B$ at a given vertex $t_0 \in \TT_n$. 
\end{corollary} 
 
In the geometric type case, there is a distinguished choice 
of ground ring for the corresponding cluster algebra: take $\AAA$ to be 
the polynomial ring $\ZZ[p_i : i \in I']$. 
In the situation of Corollary~\ref{cor:geometric type local}, 
we will denote the corresponding cluster algebra $\AA_\AAA (\MM)$ simply by 
$\AA (\tilde B)$. 
In the notation of Section~\ref{sec:setup}, $\AA (\tilde B)$ is the subring of 
the ambient field $\FFcal$ generated by cluster variables $x_j (t)$ 
for all $j \in I$ and $t \in \TT_n$ together with the 
generators $p_i \, (i \in I')$ of $\PP$. 
 
The set $I'$ is allowed to be empty: this simply means that
all the coefficients $p_j (t)$ in the corresponding exchange pattern are equal to~$1$. 
In this case, we have $\tilde B(t) = B(t)$ for all~$t$.
 
 We note that the class of exchange patterns of geometric type 
(and the corresponding cluster algebras $\AA (\tilde B)$) 
is stable under the operations of restriction and direct product 
introduced in Section~\ref{sec:setup}. 
The restriction from $I$ to a subset $J$ amounts to removing 
from $\tilde B$ the columns labeled by $I - J$; 
the direct product operation replaces two matrices 
$\tilde B_1$ and $\tilde B_2$ by the matrix 
$\bmat{\tilde B_1}{0}{0}{\tilde B_2}$.

\section{The rank~$2$ case} 
\label{sec:rank 2} 
 
In this section, we illustrate 
%Theorems~\ref{th:laurent-binomial} and \ref{th:cartan-killing} 
the above results and constructions by treating in detail the special case $n = 2$. 
We label vertices and edges of the tree $\TT_2$ as in 
(\ref{eq:tree-2}). 
For $m\in\ZZ$, it will be convenient to denote by $\rem{m}$ 
the element of $\{1,2\}$ congruent to $m$ modulo~$2$. 
Thus, $\TT_2$ consists of vertices $t_m$ and edges 
$t_m \overunder{\rem{m}}{} t_{m+1}$ for all $m \in \ZZ$. 
We use the notation of Example~\ref{ex:rank 2}, so the clusters 
are of the form ${\bf x} (t_m) = \{y_m, y_{m+1}\}$, and the exchange 
relations are given by (\ref{eq:exchange rank 2}), with coefficients 
$q_m$ and $r_m$ satisfying (\ref{eq:qqr=rr}). 
More specifically, we have 
\begin{equation} 
\label{eq:rank 2 nomenclature} 
x_{\rem{m}} (t_m) = y_m, \,\, x_{\rem{m+1}} (t_m) = y_{m+1}, \,\, 
p_{\rem{m}} (t_m) = q_{m+1}, \,\, 
p_{\rem{m+1}} (t_m) = r_{m} 
\end{equation} 
for $m \in \ZZ$. 
Choose ${\bf x}(t_1) = \{y_1, y_{2}\}$ as the \emph{initial} cluster. 
According to Theorem~\ref{th:laurent-binomial}, each cluster variable $y_m$ 
can be expressed as a Laurent polynomial in $y_1$ and $y_2$ with 
coefficients in $\ZZ \PP$. 
Let us write this Laurent polynomial as 
\begin{equation} 
\label{eq:denominator} 
y_m 
= \frac{N_m(y_1,y_2)}{y_1^{d_1(m)}y_2^{d_2(m)}}\,, 
\end{equation} 
where $N_m(y_1,y_2) \in \ZZ \PP [y_1, y_2]$ is a polynomial not divisible by $y_1$ or $y_2$. 
We will investigate in detail the denominators of 
these Laurent polynomials. 
 
Recall that for $n = 2$ the matrices $B(t)$ 
are given by (\ref{eq:B(t)-rank2}). 
Thus, all these matrices have the same associated Cartan matrix 
\begin{equation} 
\label{eq:cartan-rank2} 
A=A(B(t))=\bmat{2}{-b}{-c}{2} 
\end{equation} 
(see (\ref{eq:cartan matrix})). 
We will show that the denominators in (\ref{eq:denominator}) have a nice 
interpretation in terms of the root system associated to $A$. 
Let us recall some basic facts about this root system (cf., e.g., \cite{kac}). 
 
Let $Q \cong \ZZ^2$ be a lattice of rank $2$ with a fixed basis 
$\{\alpha_1, \alpha_2\}$ of \emph{simple roots}. 
The \emph{Weyl group} $W = W(A)$ is a group of linear transformations of 
$Q$ generated by two \emph{simple reflections} $s_1$ and $s_2$ 
whose action in the basis of simple roots is given by 
\begin{equation} 
\label{eq:s1-s2} 
s_1 = \bmat{-1}{b}{0}{1}\ ,\quad 
s_2 = \bmat{1}{0}{c}{-1}\ . 
\end{equation} 
Since both $s_1$ and $s_2$ are involutions, each element of $W$ 
is one of the following: 
$$w_1 (m) = s_1 s_2 s_1 \cdots s_{\rem{m}}, \,\, 
w_2 (m) = s_2 s_1 s_2 \cdots s_{\rem{m+1}} \ ;$$ 
here both products are of length $m \geq 0$. 
It is well known that $W$ is finite if and only if $bc \leq 3$; 
we shall refer to this as the \emph{finite case}. 
The \emph{Coxeter number} $h$ of $W$ is the order of $s_1 s_2$ in 
$W$; it is given by Table~\ref{tab:Coxeter number}. 
In the finite case, $W$ is the dihedral group of order $2h$, and 
its elements can be listed as follows: $w_1 (0) = w_2 (0) = e$ 
(the identity element), $w_1 (h) = w_2 (h) = w_0$ (the longest 
element), and $2h - 2$ distinct elements $w_1 (m), w_2 (m)$ 
for $0 < m < h$. 
In the infinite case, all elements $w_1 (m)$ and $w_2 (m)$ for $m > 0$
are distinct. 
%\smallskip 
\begin{table}[h] 
\begin{center} 
\begin{tabular}{|c|c|c|c|c|c|} 
\hline 
$bc$ & 0 & 1 & 2 & 3 & $\geq 4$ \\ 
\hline 
$h$ & 2 & 3 & 4 & 6 & $\infty$ \\ 
\hline 
\end{tabular} \,. 
\end{center} 
\medskip 
\caption{The Coxeter number} 
\label{tab:Coxeter number} 
\end{table} 
%\smallskip 

A vector $\alpha \in Q$ is a \emph{real root} for $A$ if 
it is $W$-conjugate to a simple root. 
Let $\Phi$ denote the set of real roots for $A$. 
It is known that 
$\Phi = \Phi_+ \cup (- \Phi_+)$, 
where 
$$\Phi_+ = 
\{\alpha = d_1 \alpha_1 + d_2 \alpha_2 \in \Phi : \ d_1, d_2 \geq 0\}$$ 
is the set of positive real roots. 
In the finite case, $\Phi_+$ has cardinality $h$, 
and we have 
$$\Phi_+ = \{w_1 (m) \alpha_{\rem{m+1}}: 0 \leq m < h\} 
%= \{w_2 (m) \alpha_{\rem{m+2}}: 0 \leq m < h\} 
\ .$$ 
In the infinite case, we have 
$$\Phi_+ = \{w_1 (m) \alpha_{\rem{m+1}}, \,\, 
w_2 (m) \alpha_{\rem{m+2}}: m \geq 0\} \ ,$$ 
with all the roots $w_1 (m) \alpha_{\rem{m+1}}$ and 
$w_2 (m) \alpha_{\rem{m+2}}$ distinct.

We will represent the denominators in (\ref{eq:denominator}) 
as vectors in the root lattice $Q$ by setting 
$$\delta (m) = d_1 (m) \alpha_1 + d_2 (m) \alpha_2$$ 
for all $m \in \ZZ$. 
In particular, we have $\delta (1) = - \alpha_1$ and 
$\delta (2) = - \alpha_2$. 
 
\begin{theorem} 
\label{th:denominators-rank 2} 
In the rank $2$ case (finite or infinite alike), 
cluster variables are uniquely up to a multiple from $\PP$ 
determined by their denominators in the Laurent expansions with 
respect to a given cluster. 
The set of these denominators is naturally identified with 
$\{- \alpha_1, - \alpha_2\} \cup \Phi_+$. 
More precisely:\\
{\rm (i)} In the infinite case, we have 
\begin{equation} 
\label{eq:denominator-roots-infinite} 
\delta (m+3) =  w_1 (m) \alpha_{\rem{m+1}}, \quad 
\delta (-m) =  w_2 (m) \alpha_{\rem{m+2}} \quad \quad 
(m \geq 0) \ . 
\end{equation} 
In particular, all $y_m$ for $m \in \ZZ$ have different denominators 
$y_1^{d_1(m)} y_2^{d_2(m)}$.\\ 
{\rm (ii)} In the finite case, we have 
\begin{equation} 
\label{eq:denominator-roots-finite} 
\delta (m+3) =  w_1 (m) \alpha_{\rem{m+1}} \quad \quad 
(h > m \geq 0) \ , 
\end{equation} 
and $\delta (m + h + 2) = \delta (m)$ for all $m \in \ZZ$, so the 
denominators $y_1^{d_1(m)} y_2^{d_2(m)}$ are periodic with the period $h + 2$. 
Moreover, $y_{m+h+2}/y_m \in \PP$ for $m \in \ZZ$. 
\end{theorem} 
 
\proof 
Before giving a unified proof of 
Theorem~\ref{th:denominators-rank 2}, 
we notice that part (ii) can be proved by a direct calculation. 
We will express all coefficients in terms of $r_1$, $r_2$ and $q_m$ for $m \in \ZZ$ 
(cf. Proposition~\ref{pr:free generators universal}). 
With the help of \texttt{Maple} we find that in each finite case, 
the elements $y_m$ for $3 \leq m \leq h+4$ are the following 
Laurent polynomials in $y_1$ and $y_2$. 
 
\smallskip 
 
\noindent {\bf Type $A_1 \times A_1$: $b = c = 0$.} 
$$y_3 = \frac{q_2+r_2}{y_1}, \quad y_4 = \frac{q_3(q_1+r_1)}{r_1 y_2}, \quad 
y_5 = \frac{q_4 y_1}{r_2}, \quad 
y_6 = \frac{r_1 q_5 y_2}{q_1 q_3} \ .$$ 
 
\smallskip 
 
\noindent {\bf Type $A_2$: $b = c = 1$.} 
$$y_3 = \frac{q_2+r_2}{y_1}, \quad y_4 = 
\frac{q_3 (r_2 (q_1 y_1 + r_1) + r_1 q_2 y_2)}{r_1 y_1 y_2},$$ 
$$y_5 = \frac{q_3 q_4 (q_1 y_1 + r_1)}{r_1 y_2}, \quad 
y_6 = \frac{q_4 q_5 y_1}{r_2}, \quad 
y_7 = \frac{r_1 q_5 q_6 y_2}{r_2 q_1 q_3} \ .$$ 
 
\smallskip 
 
\noindent {\bf Type $B_2$: $b = 1, c = 2$.} 
$$y_3 = \frac{q_2 y_2^2 + r_2}{y_1}, \quad 
y_4 = \frac{q_3 (r_2 (q_1 y_1 + r_1) + r_1 q_2 y_2^2)}{ r_1 y_1 y_2},$$ 
$$y_5 = \frac{q_3^2 q_4 (r_2 (q_1 y_1 + r_1)^2 + 
r_1^2 q_2 y_2^2)}{r_1^2 y_1 y_2^2}, \quad 
y_6 = \frac{q_3 q_4 q_5 (q_1 y_1 + r_1)}{r_1 y_2},$$ 
$$y_7 = \frac{q_4 q_5^2 q_6 y_1}{r_2}, \quad 
y_8 = \frac{r_1 q_5 q_6 q_7 y_2}{r_2 q_1 q_3} \ .$$ 
 
\smallskip 
 
\noindent {\bf Type $C_2$: $b = 2, c = 1$.} 
$$y_3 = \frac{q_2 y_2 + r_2}{y_1}, \quad 
y_4 = \frac{q_3 (r_2^2 q_1 y_1^2 + r_1 (q_2 y_2 + r_2)^2)}{r_1 y_1^2 y_2},$$ 
$$y_5 = \frac{q_3 q_4 (r_2 q_1 y_1^2 + r_1(q_2 y_2 + r_2))}{r_1 y_1 y_2}, 
\quad y_6 = \frac{q_3 q_4^2 q_5 (q_1 y_1^2 + r_1)}{r_1 y_2},$$ 
$$y_7 = \frac{q_4 q_5 q_6 y_1}{r_2}, \quad 
y_8 = \frac{r_1 q_5 q_6^2 q_7 y_2}{r_2^2 q_1 q_3} \ .$$ 
 
\smallskip 
 
\noindent {\bf Type $G_2$: $b = 1, c = 3$.} 
$$y_3 = \frac{q_2 
y_2^3 + r_2}{y_1}, \quad y_4 = \frac{q_3 (r_2 (q_1 y_1 + r_1) + 
r_1 q_2 y_2^3)}{ r_1 y_1 y_2},$$ 
$$y_5 = \frac{q_3^3 q_4 (r_2^2 (q_1 y_1 + r_1)^3 + 
r_1^2 q_2 y_2^3 (r_1 q_2 y_2^3 + 3 r_2 q_1 y_1 + 2 r_1 r_2))} 
{r_1^3 y_1^2 y_2^3},$$ 
$$y_6 = \frac{q_3^2 q_4 q_5 (r_2 (q_1 y_1 + r_1)^2 + r_1^2 q_2 
y_2^3)}{r_1^2 y_1 y_2^2}, \quad y_7 = \frac{q_3^3 q_4^2 q_5^3 q_6 
(r_2 (q_1 y_1 + r_1)^3 + r_1^3 q_2 y_2^3)}{r_1^3 y_1 y_2^3},$$ 
$$y_8 = \frac{q_3 q_4 q_5^2 q_6 q_7 (q_1 y_1 + r_1)}{r_1 y_2}, 
\quad y_9 = \frac{q_4 q_5^3 q_6^2 q_7^3 q_8 y_1}{r_2}, 
\quad y_{10} = \frac{r_1 q_5 q_6 q_7^2 q_8 q_9 y_2}{r_2 q_1 q_3} 
\ .$$ 
 
\smallskip 
 
\noindent {\bf Type $G_2^\vee$: $b = 3, c = 1$.} 
$$y_3 = \frac{q_2 y_2 + r_2}{y_1}, \quad 
y_4 = \frac{q_3 (r_1 (q_2 y_2 + r_2)^3 + 
r_2^3 q_1 y_1^3)}{ r_1 y_1^3 y_2} \ ,$$ 
$$y_5 = \frac{q_3 q_4 (r_1 (q_2 y_2 + r_2)^2 + 
r_2^2 q_1 y_1^3)}{ r_1 y_1^2 y_2} \ ,$$ 
$$y_6 = \frac{q_3^2 q_4^3 q_5 (r_1^2 (q_2 y_2 + r_2)^3 + 
r_2^2 q_1 y_1^3 (r_2 q_1 y_1^3 + 3 r_1 q_2 y_2 + 2 r_1 r_2))} 
{r_1^2 y_1^3 y_2^2},$$ 
$$y_7 = \frac{q_3 q_4^2 q_5 q_6 
(r_2 (q_1 y_1^3 + r_1) + r_1 q_2 y_2)}{r_1 y_1 y_2}, \quad 
y_8 = \frac{q_3 q_4^3 q_5^2 q_6^3 q_7 (q_1 y_1^3 + r_1)}{r_1 y_2} 
\ ,$$ 
$$y_9 = \frac{q_4 q_5 q_6^2 q_7 q_8 y_1}{r_2}, 
\quad y_{10} = \frac{r_1 q_5 q_6^3 q_7^2 q_8^3 q_9 y_2}{r_2^3 q_1 q_3} 
\ .$$ 
%\smallskip 
By inspection, the exponents of the $y$ variables in the denominators above agree 
with (\ref{eq:denominator-roots-finite}). 
Furthermore, we see that $y_{m+h+2}/y_m \in \PP$ for $m = 1,2$, 
and therefore for any $m$. 
 
Now let us give a unified proof of both parts of 
Theorem~\ref{th:denominators-rank 2}. 
For $\delta = d_1 \alpha_1 + d_2 \alpha_2 \in Q$, let us denote 
\[ 
[\delta]_+ = \max (d_1, 0) \alpha_1 + \max (d_2, 0) \alpha_2. 
\] 
It then follows from the relations (\ref{eq:exchange rank 2}) that 
\begin{equation} 
\label{eq:delta-recurrence} 
\delta(m+1) + \delta(m-1) = 
\begin{cases} 
b [\delta (m)]_+ & \text{if $m$ is odd;} \\[.15in] 
c [\delta (m)]_+ & \text{if $m$ is even.} 
\end{cases} 
\end{equation} 
Starting with $\delta (1) = - \alpha_1$ and 
$\delta (2) = - \alpha_2$, we use (\ref{eq:delta-recurrence}) to 
see that 
$$\delta (3) = \alpha_1, \quad \delta (4) = b \alpha_1 + \alpha_2 
= s_1 \alpha_2 \ ,$$ 
and also 
$$\delta (0) = \alpha_2, \quad \delta (-1) = \alpha_1 + c \alpha_2 
= s_2 \alpha_1 \ .$$ 
This proves (\ref{eq:denominator-roots-infinite}) and 
(\ref{eq:denominator-roots-finite}) for $m = 0$ and $m = 1$. 
 
Proceeding by induction on $m$, let us now assume that $m \geq 2$, 
and that both vectors $\delta (m+1)$ and $\delta (m+2)$ are 
positive roots given by %(\ref{eq:denominator-roots-infinite}) and 
(\ref{eq:denominator-roots-finite}); in the finite case, we also 
assume that $m \leq h$. 
If $m$ is odd, it follows from (\ref{eq:delta-recurrence}) that 
\begin{eqnarray*} 
&&\delta (m+3) = b \delta (m+2) - \delta (m+1) = 
b w_1 (m-1) \alpha_1 - w_1 (m-2) \alpha_2 \\ 
&&=w_1 (m-1) (s_1 \alpha_2 - \alpha_2) - w_1 (m-2) \alpha_2 \\ 
&&= w_1 (m) \alpha_2 - w_1 (m-2) (s_2 \alpha_2 + \alpha_2) = 
w_1 (m) \alpha_2 \ . 
\end{eqnarray*} 
Thus, $\delta (m+3)$ is also given by 
(\ref{eq:denominator-roots-finite}). 
The same argument works with $m$ even, and, in the infinite 
case, with $m$ negative. 
This completes the proof of (\ref{eq:denominator-roots-infinite}) 
and~(\ref{eq:denominator-roots-finite}). 
 
We conclude by a unified proof of the periodicity $\delta (m + h + 2) = \delta (m)$ 
in the finite case. 
Without loss of generality, we assume that $m = 1$. 
By the above inductive argument, we know that 
$\delta (h + 3)$ is given by (\ref{eq:denominator-roots-finite}). 
Thus, we have $\delta (h + 3) =  w_1 (h) \alpha_{\rem{h+1}} 
= w_0 \alpha_{\rem{h+1}}$. 
It is known---and easy to check---that 
$w_0 = \bmat{0}{-1}{-1}{0}$ if $b = c = 1$ (the only case when $h$ is odd), and 
$w_0 = \bmat{-1}{0}{0}{-1}$ in the other three finite cases. 
Hence $w_0 \alpha_{\rem{h+1}} = - \alpha_1$ in all cases, and we are done. 
\endproof

In the rest of the section, we discuss the normalized case 
(see Definition~\ref{def:normalized}). 
In this context, the periodicity property in 
Theorem~\ref{th:denominators-rank 2} (ii) can be sharpened 
considerably. 
 
\begin{proposition} 
\label{pr:periodicity rank 2 normalized} 
In the normalized finite case of a cluster algebra of rank $2$, 
we have $q_{m+h + 2} = q_m$, $r_{m+h + 2} = r_m$, and 
$y_{m+h + 2} = y_m$ for all $m \in \ZZ$. 
\end{proposition} 
 
\proof 
We set $u_m = q_m/r_m$ for $m \in \ZZ$. 
Then $q_m$ and $r_m$ are recovered from $u_m$~by 
\begin{equation} 
\label{eq:qru} 
q_m = \frac{u_m}{1 \oplus u_m}, \quad r_m = \frac{1}{1 \oplus u_m} 
\end{equation}  
(cf.~(\ref{eq:p through u})). 
To establish the periodicity of $q_m$ and $r_m$, it suffices 
to show that $u_{m+h + 2} = u_m$. 
Furthermore, it is enough to 
check that $u_{h+3} = u_1$ and $u_{h+4} = u_2$ in each of the 
following four cases: $(b,c) = (0,0)$ (type $A_1 \times A_1$), $(1,1)$ 
(type $A_2$), $(1,2)$ (type $B_2$), and $(1,3)$ (type $G_2$). 
This is done by direct calculation using the relation 
\begin{equation} 
\label{eq:u-recurrence} 
u_{m+1} u_{m-1} = 
\begin{cases} 
(1 \oplus u_m)^c & \text{if $m$ is odd;} \\[.15in] 
(1 \oplus u_m)^b & \text{if $m$ is even,} 
\end{cases} 
\end{equation} 
which, in view of (\ref{eq:rank 2 nomenclature}), is a consequence of (\ref{eq:uu'}) and 
(\ref{eq:3 relation u}). 
Below are the intermediate steps of this calculation, which in each case expresses 
$u_m$ for $3 \leq m \leq h+2$ as a rational function 
in $u_1$ and $u_2$, and then confirms $u_{h+3} = u_1$ and $u_{h+4} = u_2$.

\smallskip 
 
\noindent {\bf Type $A_1 \times A_1$: $b = c = 0$, $h = 2$.} 
$$u_3 = \frac{1}{u_1}, \quad u_4 = \frac{1}{u_2} 
%, \quad u_5 = u_1, \quad u_6 = u_2 
\ .$$ 
 
\smallskip 
 
\noindent {\bf Type $A_2$: $b = c = 1$, $h = 3$.} 
$$u_3 = \frac{1 \oplus u_2}{u_1}, \quad u_4 = 
\frac{1 \oplus u_1 \oplus u_2}{u_1 u_2}, \quad u_5 = \frac{1 \oplus u_1}{u_2} 
%, \quad u_6 = u_1, \quad u_7 = u_2 
\ .$$ 
 
\smallskip 
 
\noindent {\bf Type $B_2$: $b = 1, c = 2$, $h = 4$.} 
$$u_3 = \frac{1 \oplus u_2}{u_1}, \quad u_4 = 
\frac{(1 \oplus u_1 \oplus u_2)^2}{u_1^2 u_2}, \quad 
u_5 = \frac{(1 \oplus u_1)^2 \oplus u_2}{u_1 u_2}, \quad 
u_6 =  \frac{(1 \oplus u_1)^2}{u_2} \ .$$

\smallskip

\noindent {\bf Type $G_2$: $b = 1, c = 3$, $h = 6$.} 
$$u_3 = \frac{1 \oplus u_2}{u_1}, \quad u_4 = 
\frac{(1 \oplus u_1 \oplus u_2)^3}{u_1^3 u_2}, \quad 
u_5 = \frac{(1 \oplus u_1)^3 \oplus u_2(3 u_1 \oplus u_2 \oplus 2)}{u_1^2 
u_2}\ ,$$ 
$$u_6 =  \frac{((1 \oplus u_1)^2 \oplus u_2)^3}{u_1^3 u_2^2}, 
\quad u_7 =  \frac{(1 \oplus u_1)^3 \oplus u_2}{u_1 u_2}, 
\quad u_8 =  \frac{(1 \oplus u_1)^3}{u_2} \ .$$

\smallskip 
 
To complete the proof of Proposition~\ref{pr:periodicity rank 2 normalized}, 
it remains to show that $y_{h+3} = y_1$ and $y_{h+4} = y_2$ in each of the four 
cases. 
By Theorem~\ref{th:denominators-rank 2} (ii), both ratios 
$y_{h+3}/y_1$ and $y_{h+4}/y_2$ belong to $\PP$, and we only need 
to show that these two elements of $\PP$ are equal to $1$ in the 
normalized case. 
The ratios in question were explicitly computed in the course of 
the proof of Theorem~\ref{th:denominators-rank 2}. 
We see that it all boils down to the following identities: 
\begin{eqnarray}
\label{eq:r-via-q-A1A1}
&&\text{{\rm Type $A_1 \times A_1$:} $r_1 = q_3$.}\\ 
\label{eq:r-via-q-A2}
&&\text{{\rm Type $A_2$:} $r_1 = q_3 q_4$.}\\
\label{eq:r-via-q-B2}
&&\text{{\rm Type $B_2$:} $r_1 = q_3 q_4 q_5$, $r_2 = q_4 q_5^2 q_6$.}\\ 
\label{eq:r-via-q-G2}
&&\text{{\rm Type $G_2$:} $r_1 = q_3 q_4 q_5^2 q_6 q_7$, 
$r_2 = q_4 q_5^3 q_6^2 q_7^3 q_8$.}
\end{eqnarray} 
 
All these identities can be proved by a direct calculation 
(preferably, with a computer): first 
express both sides in terms of the $u_m$ using (\ref{eq:qru}), 
and then replace each $u_m$ by 
the above expression in terms of $u_1$ and $u_2$. 
\endproof 
 
\begin{remark} 
\label{rem:Zamolodchikov rank 2} 
{\rm Periodicity of the recurrence (\ref{eq:u-recurrence}) 
is a very special case of the periodicity phenomenon for $Y$-system 
recurrences in the theory of the thermodynamic Bethe ansatz~\cite{zamolodchikov}. 
We plan to address the case of arbitrary rank in a forthcoming paper.} 
\end{remark} 
 
%\begin{remark} 
%\label{rem:dual roots rank 2} 
%{\rm 
We have shown that in a finite normalized case, 
any coefficient $r_m$ can be written as a monomial in 
$q_{m+2}, \dots, q_{m + h}$. 
There is a nice uniform way to write down these monomials 
using the \emph{dual root system} of the system $\Phi$ considered 
above. 
Recall that the dual root system $\Phi^\vee$ is %the root system 
associated with the transposed Cartan matrix. 
The simple roots of $\Phi^\vee$ are called \emph{simple coroots} and denoted 
$\alpha_1^\vee$ and $\alpha_2^\vee$. 
They generate the \emph{coroot lattice} $Q^\vee$. 
The Weyl groups of $\Phi$ and $\Phi^\vee$ are naturally isomorphic 
to each other, and it is common to identify them. 
The same Weyl group $W$ acts on $Q^\vee$ 
so that the action of simple reflections $s_1$ and $s_2$ 
in the basis of simple coroots is given by the matrices 
$$s_1 = \bmat{-1}{c}{0}{1}\ ,\quad 
s_2 = \bmat{1}{0}{b}{-1}$$ 
(cf. (\ref{eq:s1-s2})). 
By inspection (cf.\ (\ref{eq:r-via-q-A1A1})--\ref{eq:r-via-q-G2})), 
one obtains the following unified expression for $r_m$. 

\begin{proposition}
\label{pr:dual roots rank 2}
For $i \in \{1,2\}$ and $m \geq 0$, let $c(i,m)$ denote the 
coefficient of $\alpha_i^\vee$ in the root $w_i (m) \alpha_{\rem{i+m}}^\vee$. 
Then, in every finite case, we have 
\begin{equation} 
\label{eq:r through q} 
r_k = \prod_{m = 0}^{h-2} q_{k + m + 2}^{c(\rem{k+1},m)} 
\end{equation} 
for all $k \in \ZZ$. 
\end{proposition} 
 
\begin{remark} 
\label{rem:normalized coefficient group} 
{\rm The relations (\ref{eq:r through q}) together with the 
periodicity relation $q_{m + h + 2} = q_m$ imply that in each 
finite normalized case the subgroup of the coefficient group 
$\PP$ generated by all $q_m$ and $r_m$ is in fact generated  
by $q_1, \dots, q_{h+2}$. 
Thus, this group is dramatically different from the universal 
coefficient group. 
One might ask if the normalization condition imposes some 
relations among $q_1, \dots, q_{h+2}$. 
The answer to this question is negative: it is possible 
to take $\PP=\Trop(q_1,\dots,q_{h+2})$ 
(cf.\ Example~\ref{ex:coefficients-tropical}), so that 
multiplicatively, $\PP$ is a free abelian group generated 
by $q_1, \dots, q_{h+2}$. 
%To see this, make this multiplicative group $\PP$ into a semifield 
%as in Example~\ref{ex:coefficients-tropical}, with respect to the generators 
%$q_1, \dots, q_{h+2}$. 
With the help of (\ref{eq:r through q}) and the 
periodicity relation $q_{m + h + 2} = q_m$, we can view all 
$r_m$  and $q_m$ for $m \in \ZZ$ as elements of $\PP$. 
By the definition of auxiliary addition (see (\ref{eq:tropical addition})), we have 
$r_m \oplus q_m = 1$ for all~$m$. 
Hence these coefficients give 
rise to a normalized exchange pattern (which is in fact 
of geometric type according to Definition~\ref{def:geometric type}). 
 
In the infinite case, the situation is very different: 
there exists a normalized exchange pattern such that the 
coefficients $r_m$ and $q_m$ generate the universal coefficient 
group (see Proposition~\ref{pr:free generators universal}). 
To see this, take $\PP$ to be the group of formal infinite 
products $\prod_{m \in \ZZ} q_m^{c_m}$, where the exponents 
$c_m$ can be arbitrary integers. 
Make $\PP$ into a semifield by setting 
$$\prod_{m \in \ZZ} q_m^{c_m} \oplus \prod_{m \in \ZZ} q_m^{c'_m} 
= \prod_{m \in \ZZ} q_m^{\min (c_m, c'_m)} \ .$$ 
Taking (\ref{eq:r through q}) as an inspiration, 
let us define the elements $r_k \in \PP$ for all $k \in \ZZ$ by 
$$r_k = \prod_{m \geq 0} (q_{k + m + 2} q_{k - m - 2})^{c(\rem{k+1},m)} 
\ ,$$ 
with the same exponents as in (\ref{eq:r through q}). 
Also, view each $q_k$ as an element of $\PP$ via 
$q_k = \prod_{m \in \ZZ} q_m^{\delta_{km}}$. 
A direct check then shows that these 
coefficients give rise to a normalized exchange pattern 
(that is, satisfy (\ref{eq:qqr=rr})), and also that 
in $\PP$ there are no relations among the elements $r_1$, $r_2$ 
and $q_k\,$, for $k \in \ZZ$.} 
\end{remark}

\section{The exchange graph} 
\label{sec:cluster graph} 
 
In this section, we restrict our attention to \emph{normalized} exchange 
patterns $\MM = (M_j (t) :  j \in I, t \in \TT_n)$. 
Let $(x_j (t)  :  j \in I, t \in \TT_n)$ be the 
family of cluster variables associated to $\MM$. 
As in Section~\ref{sec:setup}, we view them 
as elements of the ambient field~$\FFcal$.

A normalized exchange pattern gives rise to a natural equivalence relation on 
the set of vertices of $\TT_n$. 
Informally speaking, two vertices are equivalent with respect to 
$\MM$ if they have the same clusters and the same exchange 
relations, up to relabeling  of cluster variables. 
Here is a precise definition. 
 
\begin{definition} 
\label{def: cluster equivalence} 
{\rm We say that two vertices $t$ and $t'$ in $\TT_n$ 
are $\MM$-equivalent if there exists a permutation 
$\sigma$ of the index set $I$ such that: 
 
\smallskip 
 
\noindent (1) $x_i (t') = x_{\sigma (i)} (t)$ for all $i \in I$; 
 
\smallskip 
 
\noindent (2) if $t \overunder{\sigma (j)}{} t_1$ and 
$t' \overunder{j}{} t'_1$ for some $j \in I$, then 
$M_j (t')({\bf x} (t')) = M_{\sigma (j)} (t)({\bf x}(t))$ 
and $M_j (t'_1)({\bf x} (t')) = M_{\sigma (j)} (t_1)({\bf x}(t))$. 
%$u_j (t') = u_{\sigma (j)} (t)$ for all $j \in I$. 
} 
\end{definition} 

\begin{remark}
\label{rem:exchange-uniqueness}
{\rm We believe that condition (2) in Definition~\ref{def: cluster equivalence} 
is in fact a consequence of (1).
Thus, once a cluster has been fixed as a \emph{set}
(forgetting the labels and the cluster's location in the tree $\TT_n$), 
the exchange relations involving its elements are uniquely determined.  
}
\end{remark}

In view of Proposition~\ref{pr:normalized}, a normalized exchange pattern 
can be defined in terms of a family of integer matrices 
$B(t)$ together with a family of elements $u_j (t) \in \PP$ 
satisfying relations (\ref{eq:uu'})--(\ref{eq:3 relation u}). 
In these terms, property (2) above can be rephrased as saying that 
$b_{ij} (t') = b_{\sigma (i), \sigma (j)} (t)$ and 
$u_j (t') = u_{\sigma (j)} (t)$ for all $i, j \in I$. 
 
The following property is an immediate consequence of 
Proposition~\ref{pr:normalized}.
 
\begin{proposition} 
\label{pr:equivalence adjacency} 
Suppose $t$ and $t'$ are $\MM$-equivalent. 
For every vertex $t_1$ adjacent to $t$, there is a unique vertex 
$t'_1$ adjacent to $t'$ and equivalent to $t_1$. 
In the notation of Definition~{\rm \ref{def: cluster equivalence},} 
if $t \overunder{\sigma (j)}{} t_1$ for some $j \in I$, then 
$t' \overunder{j}{} t'_1$. 
\end{proposition} 
 
\begin{definition} 
\label{def:exchange graph} 
{\rm The \emph{exchange graph} $\Gamma_\MM$ associated with a normalized 
exchange pattern $\MM$ has the $\MM$-equivalence classes as 
vertices, joined by an edge if they have 
adjacent representatives in $\TT_n$.} 
\end{definition} 
 
As an immediate corollary of Proposition~\ref{pr:equivalence adjacency}, 
we obtain the following. 
 
\begin{corollary} 
\label{cor:exchange graph regular} 
The exchange graph is $n$-regular, i.e., every vertex has 
precisely $n$ edges emanating from it. 
\end{corollary} 
 
Note that passing from $\TT_n$ to $\Gamma_\MM$ may result in 
losing the coloring of edges by elements of the index set $I$. 
This will happen if a permutation $\sigma$ in 
Definition~\ref{def: cluster equivalence} is non-trivial. 
(This permutation can be thought of as ``discrete monodromy" 
of our graph.) 

\begin{example}
\label{ex:rank 2 exchange graph}
{\rm The results in Section~\ref{sec:rank 2} 
(more specifically, Theorem~\ref{th:denominators-rank 2} 
and Proposition~\ref{pr:periodicity rank 2 normalized}) 
show that, in the rank $2$ case, the $\MM$-equivalence and the exchange
graph can be described as follows.
In the infinite case, no two distinct vertices of $\TT_2$ are 
$\MM$-equivalent to each other, so the exchange graph is $\TT_2$ itself.
In the finite case, $t_m$ and $t_{m'}$ are $\MM$-equivalent if and only if 
$m \equiv m' \bmod (h+2)$, so the exchange graph is a cycle of length~$h+2$.
For example, in the type $A_2$ case,
we have $h=3$, and the graph is a pentagon.
Figure~\ref{fig:pentagon-A2} shows the corresponding exchange graph,
together with its clusters and exchange relations,
written with the help of (\ref{eq:exchange rank 2})
and (\ref{eq:r-via-q-A2}).
The discrete monodromy is present in this special case, 
as the variables $y_1$ and $y_2$ get switched after a full 5-cycle of
exchanges:
\[
\begin{pmatrix}y_1\\y_2\end{pmatrix}\to
\begin{pmatrix}y_3\\y_2\end{pmatrix}\to
\begin{pmatrix}y_3\\y_4\end{pmatrix}\to
\begin{pmatrix}y_5\\y_4\end{pmatrix}\to
\begin{pmatrix}y_5\\y_1\end{pmatrix}\to
\begin{pmatrix}y_2\\y_1\end{pmatrix}. 
\]
%In the geometric case, we recover 
The corresponding cluster algebra can be realized 
as the homogeneous coordinate ring of the Grassmannian $Gr_{2,5}$; 
see Example~\ref{ex:rank 2}. 
%Figure~\ref{fig:pentagon}. 

\begin{figure}[ht] 
\begin{center} 
\setlength{\unitlength}{4.2pt} 
\begin{picture}(35,35)(0,-2) 
%\thicklines 
  \put(6,0){\line(1,0){20}} 
%  \put(0,19){\line(1,0){32}} 
%  \qbezier(6,0)(11,15.5)(16,31) 
%  \qbezier(26,0)(21,15.5)(16,31) 
  \qbezier(6,0)(3,9.5)(0,19) 
  \qbezier(26,0)(29,9.5)(32,19) 
%  \qbezier(0,19)(13,9.5)(26,0) 
%  \qbezier(32,19)(19,9.5)(6,0) 
  \qbezier(0,19)(8,25)(16,31) 
  \qbezier(32,19)(24,25)(16,31) 
 
  \put(6,0){\circle*{1}} 
  \put(26,0){\circle*{1}} 
  \put(0,19){\circle*{1}} 
  \put(32,19){\circle*{1}} 
  \put(16,31){\circle*{1}} 
 
\put(2,0){\makebox(0,0){$y_5,y_1$}} 
\put(30,0){\makebox(0,0){$y_4,y_5$}} 
\put(-4,19){\makebox(0,0){$y_1,y_2$}} 
\put(36,19){\makebox(0,0){$y_3,y_4$}} 
\put(16,33){\makebox(0,0){$y_2,y_3$}} 

\put(-2,26.5){\makebox(0,0){$y_1 y_3=q_2 y_2 + q_4 q_5$}} 
\put(34,26.5){\makebox(0,0){$y_2 y_4=q_3 y_3 + q_5 q_1$}} 
\put(39.5,9){\makebox(0,0){$y_3 y_5=q_4 y_4 + q_1 q_2$}} 
\put(-7.5,9){\makebox(0,0){$y_5 y_2=q_1 y_1 + q_2 q_4$}} 
\put(16,-2.5){\makebox(0,0){$y_4 y_1=q_5 y_5 + q_2 q_3$}} 
 
\end{picture} 
\end{center} 
\caption{The exchange graph for a cluster algebra of type $A_2$} 
\label{fig:pentagon-A2} 
\end{figure}
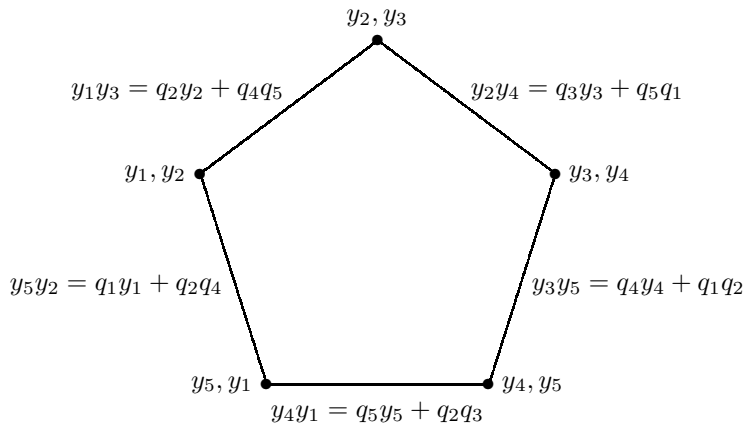 
 
We see that, in the rank $2$ case, the following conditions on 
a normalized exchange pattern $\MM$ are equivalent:
\\[.05in] 
{\rm (1)} The exchange graph $\Gamma_\MM$ is finite. 
\\[.05in] 
{\rm (2)} A cluster algebra associated to $\MM$ has 
finitely many distinct cluster variables. 
\\[.05in] 
{\rm (3)} There exists a vertex $t \in \TT_2$ 
such that the Cartan counterpart $A(B(t))$ of the 
matrix $B(t)$ %(see (\ref{eq:cartan matrix})) 
is a Cartan matrix of 
finite type, i.e., we have $|b_{12} (t) b_{21} (t)| \leq 3$.

In the sequel to this paper, we shall extend this result to 
%the case of an 
arbitrary rank. 
}
\end{example}

Returning to the general case, we shall now provide some important instances of 
$\MM$-equivalence. 
 
\begin{theorem} 
\label{th:rank 2 periodicity} 
Suppose $|b_{ij} (t) b_{ji} (t)| \leq 3$ for some vertex 
$t \in \TT_n$ and two distinct indices $i$ and $j$. 
Let $h$ be the Coxeter number of the corresponding rank $2$ root 
system, and let $t' \in \TT_n$ 
be a vertex joined with $t$ by a path of length $h+2$ whose edge 
labels alternate between $i$ and $j$. 
Then $t'$ is $\MM$-equivalent to $t$, so that the path becomes a 
cycle of length $h+2$ in the exchange graph. 
\end{theorem} 
 
\proof 
Using the restriction operation on exchange patterns, we can assume 
without loss of generality that $n = 3$. 
Let $I = \{1,2,3\}$, $i = 1$, $j = 2$, and suppose that the path
joining $t$ and $t'$ has the form 
$$t = t_1 \overunder{1}{} 
t_2 \overunder{2}{} 
t_3 \overunder{1}{} \cdots \overunder{\rem{h+2}}{} t_{h+3} = t'$$ 
(cf. (\ref{eq:tree-2})), where  $\rem{m}$ as before stands for the element 
of $\{1,2\}$ congruent to $m$ modulo $2$. 
We abbreviate $b_1 = - b_{12} (t_1)$ and $b_2 = b_{21} (t_1)$, and assume 
(without loss of generality) that $b_1$ and $b_2$ are either both positive 
or both equal to $0$. 
This agrees with the convention (\ref{eq:B(t)-rank2}) 
(also used in Section~\ref{sec:rank 2}), where 
$b_1$ was denoted by $b$, and $b_2$ by $c$. 
Recall that the value of $h$ is given by Table~\ref{tab:Coxeter number} 
in Section~\ref{sec:rank 2}; 
in particular, $h$ is even  unless $b_1 = b_2 = 1$ (the type $A_2$). 
We claim that the conditions (1) and (2) in 
Definition~\ref{def: cluster equivalence} 
hold with the permutation $\sigma$ equal to identity for $h$ even, and equal to 
the transposition of indices $1$ and $2$ in the only case when $h$ is odd. 
 
In view of the exchange property (\ref{eq:exchange1}), the cluster variable 
$x_3 (t_m)$ does not depend on $m$.  
Therefore, verifying condition (1) amounts to checking that  
$x_i (t') = x_{\sigma (i)} (t)$ for $i \in \{1,2\}$. 
Again using restriction, we can assume for the purpose of this checking 
only that $n = 2$ and $I = \{1,2\}$, in which case the required property was established in 
Proposition~\ref{pr:periodicity rank 2 normalized}.  
 
The same argument shows that the only part of condition (2) that does not follow from 
Proposition~\ref{pr:periodicity rank 2 normalized} is the equalities 
$u_3 (t') = u_3 (t)$, and  $b_{i3} (t') = b_{\sigma (i), 3} (t)$ 
for $i \in \{1,2\}$. 
Let us abbreviate $d_m = b_{\rem{m}, 3} (t_m)$ and 
$u_{m+1} = u_{\rem {m}} (t_m)$; the latter notation is chosen to be consistent 
with the notation in the proof of 
Proposition~\ref{pr:periodicity rank 2 normalized}. 
Iterating (\ref{eq:3 relation u}) (with $i$ = 3), we obtain 
\begin{equation} 
\label{eq:u-product} 
\frac{u_3 (t')}{u_3 (t)} = \prod_{m = 1}^{h+2} u_{m+1}^{\max \ (d_m,0)} 
(1 \oplus u_{m+1})^{- d_m} \ . 
\end{equation} 
Thus, we need to show that the product on the right-hand side of 
(\ref{eq:u-product}) is equal to~$1$. 
To do this, we need some preparation. 
 
Recall that the coefficients $u_m$ satisfy the relation (\ref{eq:u-recurrence}), 
which in our present notation can be rewritten as 
\begin{equation} 
\label{eq:u-chain} 
u_{m} u_{m+2} = (1 \oplus u_{m+1})^{b_{\rem{m}}} \ . 
\end{equation} 
On the other hand, the matrix mutation rules (\ref{eq:mutation})
(with $j = 3$, applied along the edges 
$t_m \overunder{\rem {m}}{} t_{m+1} \overunder{\rem {m+1}}{} t_{m+2}$),  
readily imply the following recurrence for the exponents $d_m$: 
\begin{equation} 
\label{eq:u-chain tropical} 
d_{m} + d_{m+2} = b_{\rem{m}} \max \ (d_{m+1},0) \ . 
\end{equation} 
The sequence $(u_m)$ is periodic with period 
$h+2$, by virtue of Proposition~\ref{pr:periodicity rank 2 normalized}. 
Furthermore, the relation (\ref{eq:u-chain tropical}) 
can be seen, somewhat surprisingly, as a special case of (\ref{eq:u-chain}), 
for the following version of the tropical semifield: take $\ZZ$ with 
the multiplication given by ordinary addition, and the auxiliary 
addition given by $a \oplus b = \max \ (a,b)$. 
This implies in particular that the sequence 
$(d_m)$ is also periodic with period $h+2$. 
We also note that the same periodicity holds for the sequence 
$(b_{\rem{m}})$: when $h$ is even, this is clear from the definition 
of $\rem{m}$; and in the only case when $h$ is odd, we have 
$b_1 = b_2 = 1$. 

To show that the right-hand side of 
(\ref{eq:u-product}) is equal to~$1$, we first treat the case 
$b_1 = b_2 = 0$. 
Then $h + 2 = 4$, and we have $u_{m+2} = u_m^{-1}$ and 
$d_{m+2} = - d_m$ for all~$m$. 
It follows that
\[
\begin{array}{l}
\displaystyle\prod_{m = 1}^{h+2} u_{m+1}^{\max \ (d_m,0)} 
(1 \oplus u_{m+1})^{- d_m}\\[.1in]
= \displaystyle\prod_{m = 1}^2 u_{m+1}^{\max \ (d_m,0)} 
(1 \oplus u_{m+1})^{- d_m} u_{m+1}^{- \max \ (- d_m,0)} 
(1 \oplus u_{m+1}^{-1})^{d_m}\\[.1in]
= \displaystyle\prod_{m = 1}^2 u_{m+1}^{d_m} (1 \oplus u_{m+1})^{- d_m} 
(1 \oplus u_{m+1}^{-1})^{d_m} = 1 \ ,
\end{array}
\]
as desired. 

In the remaining case when $b_1$ and $b_2$ are both positive,
we use (\ref{eq:u-chain}) and (\ref{eq:u-chain tropical}) as well as 
the above-mentioned periodicity, to obtain: 
\begin{eqnarray*} 
\left(\frac{u_3 (t')}{u_3 (t)}\right)^{b_1 b_2} 
&=& \prod_{m = 1}^{h+2} u_{m+1}^{b_1 b_2 \max \ (d_m,0)} 
(1 \oplus u_{m+1})^{- b_1 b_2 d_m}\\ 
&=& \prod_{m = 1}^{h+2} u_{m+1}^{b_1 b_2 \max \ (d_m,0)}  
(u_{m} u_{m+2})^{- b_{\rem{m+1}}d_m} \\ 
&=& u_1^{- b_{2} d_1} u_2^{b_1 ( b_2 \max \ (d_1,0) - d_2)} \\ 
&& \times \prod_{m = 1}^{h} 
u_{m+2}^{- b_{\rem{m+1}}(d_{m} + d_{m+2} - b_{\rem{m}} \max \ 
  (d_{m+1},0))}\\ 
&&\times \, u_{h+3}^{b_{\rem{h+2}} ( b_{\rem{h+1}} \max \ (d_{h+2},0) - d_{h+1})} 
u_{h+4}^{- b_{\rem{h+3}}d_{h+2}}\\ 
&=& u_1^{- b_{2} d_1} u_2^{b_1 d_0} u_{h+3}^{b_{\rem{h+4}} d_{h+3}} 
u_{h+4}^{- b_{\rem{h+3}}d_{h+2}} \\ 
&=& u_1^{(b_{\rem{h+4}}- b_{2}) d_1} u_2^{(b_1 - b_{\rem{h+3}}) d_0} = 
1\ . 
\end{eqnarray*} 
To conclude that $u_3 (t') = u_3 (t)$, just recall that the coefficient group 
$\PP$ is assumed to be torsion-free. 
 
To complete the proof of the theorem, it remains to show that 
$b_{i3} (t') = b_{\sigma (i), 3} (t)$ for $i \in \{1,2\}$. 
But this is equivalent to saying that the sequence 
$(d_m)$ is periodic with period $h+2$, which is already proven. 
\endproof 
 
We conjecture that the $\MM$-equivalence relation is generated by 
its instances described in Theorem~\ref{th:rank 2 periodicity}. 

\begin{example}
\label{ex:brick wall}
{\rm Let $\MM (B)$ be an exchange pattern of geometric type
associated with the skew-symmetric matrix 
\begin{equation}
\label{eq:brick wall matrix}
B = \left[\begin{array}{ccc} 
0 & 1 & 1 \\ 
-1 & 0 & 1\\ 
-1 & 1 & 0 \\ 
\end{array}\right] \ ,
\end{equation}
in accordance with Corollary~\ref{cor:geometric type local};
since $B$ is a $3 \times 3$ matrix, the corresponding cluster algebra is of rank $3$,
and all the coefficients $p_j (t)$ are equal to $1$. 
The exchange graph for this pattern is a ``two-layer brick wall" shown 
in Figure~\ref{fig:brick wall}. 
In this figure, distinct cluster variables 
are associated with regions: there are variables $y_m \, (m \in \ZZ)$
associated with bounded regions (``bricks"), and two more variables
$w$ and $z$ associated with the two unbounded regions.
The cluster at any vertex $t$ consists of the three cluster variables
associated to the three
regions adjacent to~$t$. 
The binomial exchange relations %corresponding to the edges 
are as follows: for $m \in \ZZ$, we have 
\begin{equation}
\label{eq:brick exchange}
\begin{array}{l}
w y_{2m} = y_{2m-1} + y_{2m+1} \ ,\\[.1in]
y_{2m-1} y_{2m+3} = y_{2m+1}^2 + w \ ,\\[.1in]
y_{m} y_{m+3} = y_{m+1} y_{m+2} + 1 \ ,\\[.1in]
y_{2m-2} y_{2m+2} = y_{2m}^2 + z \ ,\\[.1in]
y_{2m-1} z = y_{2m-2} + y_{2m} \ 
.
\end{array}
\end{equation}
To see all this, we first show that the graph in Figure~\ref{fig:brick wall}
is a cover of the exchange graph $\Gamma_\MM$. 
Pick an initial vertex $t_0 \in \Gamma_\MM$ with the matrix 
$B(t_0) = B$ given by (\ref{eq:brick wall matrix})
(abusing notation, we will use the same symbol 
for a vertex in $\TT_3$ and its image in $\Gamma_\MM$).  
Denote the cluster variables at $t_0$
by $y_1$, $y_2$ and $y_3$, so that their order agrees with 
that of rows and columns of $B$. 
Since every principal $2 \times 2$ submatrix of $B$ is of type $A_2$, 
Theorem~\ref{th:rank 2 periodicity} provides that $t_0$ is a common vertex 
of three $5$-cycles in~$\Gamma_\MM$.
These cycles are depicted in Figure~\ref{fig:brick wall}
as perimeters of the three bricks surrounding~$t_0$. 
The variable $y_i$ inside each brick indicates that this brick 
corresponds to the rank 2 exchange pattern obtained from $\MM$ via 
restriction from $I = [3]$ to $[3] \setminus \{i\}$;
equivalently, this means that $y_i$ appears in every cluster
on the perimeter of the brick. 
Now we move from $t_0$ to an adjacent vertex $t_1$ with the cluster 
$\{y_2, y_3, y_4\}$. % (cf.~Figure~\ref{fig:brick}). 
According to Proposition~\ref{prop:B(t)}, 
$$B(t_1) = \mu_1 (B) = 
\left[\begin{array}{ccc} 
0 & -1 & -1 \\ 
1 & 0 & 1\\ 
1 & -1 & 0 \\ 
\end{array}\right] \ .$$
This matrix differs from $B$ by a simultaneous cyclic permutation of rows and
columns. 
Therefore, we can apply the same construction to $t_1$, obtaining 
the three bricks surrounding this vertex. 
Continuing in the same way, we produce the entire graph in 
Figure~\ref{fig:brick wall}.
Since this graph is $3$-regular, we have covered all the vertices and edges
of the exchange graph.

\begin{figure}[ht] 
\setlength{\unitlength}{2pt}
\begin{center} 
\begin{picture}(140,55)(0,0)
%\thicklines
\multiput(0,10)(0,20){3}{\line(1,0){140}}
\multiput(10,30)(40,0){4}{\line(0,1){20}}
\multiput(30,10)(40,0){3}{\line(0,1){20}}

%\thinlines
\multiput(10,30)(40,0){4}{\circle*{2}}
\multiput(10,50)(40,0){4}{\circle*{2}}
\multiput(30,30)(40,0){3}{\circle*{2}}
\multiput(30,10)(40,0){3}{\circle*{2}}

\put(10,20){\makebox(0,0){$y_0$}}
\put(50,20){\makebox(0,0){$y_2$}}
\put(90,20){\makebox(0,0){$y_4$}}
\put(130,20){\makebox(0,0){$y_6$}}

\put(30,40){\makebox(0,0){$y_1$}}
\put(70,40){\makebox(0,0){$y_3$}}
\put(110,40){\makebox(0,0){$y_5$}}

\put(53,33){\makebox(0,0){$t_0$}}
\put(73,33){\makebox(0,0){$t_1$}}
\put(70,2){\makebox(0,0){$z$}}
\put(70,58){\makebox(0,0){$w$}}

\end{picture} 
\end{center} 
\caption{The two-layer brick wall} 
\label{fig:brick wall} 
\end{figure}
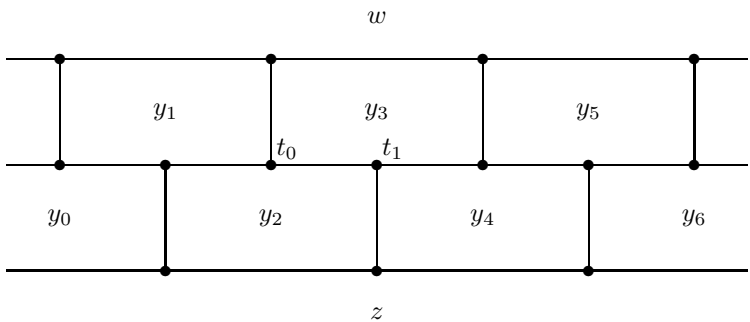 

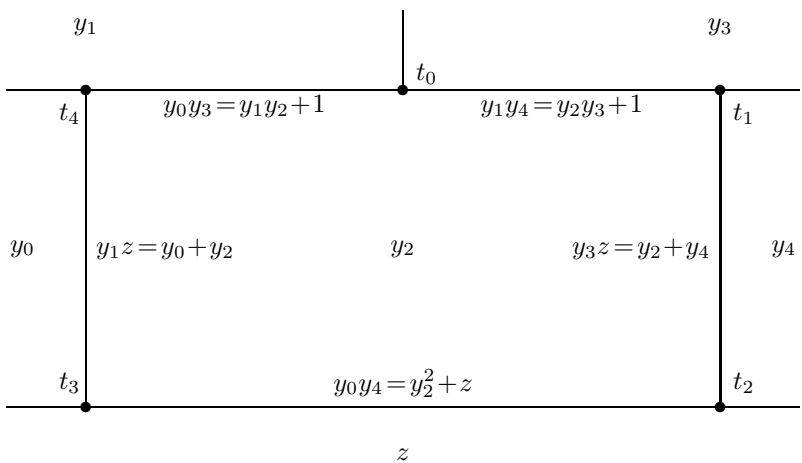
\begin{figure}[ht] 
\setlength{\unitlength}{6pt}
\begin{center} 
\begin{picture}(50,30)(0,-4)
%\thicklines
\multiput(0,0)(0,20){2}{\line(1,0){50}}
\multiput(25,20)(40,0){1}{\line(0,1){5}}
\multiput(5,0)(40,0){2}{\line(0,1){20}}

%\thinlines
\multiput(5,0)(40,0){2}{\circle*{0.7}}
\multiput(5,20)(20,0){3}{\circle*{0.7}}

\put(15,19){\makebox(0,0){$y_0 y_3\!=\!y_1 y_2 \!+\! 1$}}
\put(35,19){\makebox(0,0){$y_1 y_4\!=\!y_2 y_3 \!+\! 1$}}
\put(10,10){\makebox(0,0){$y_1 z\!=\!y_0\!+\! y_2$}}
\put(40,10){\makebox(0,0){$y_3 z\!=\!y_2\!+\! y_4$}}
\put(25,1.5){\makebox(0,0){$y_0 y_4\!=\!y_2^2\!+\! z$}}
\put(26.5,21){\makebox(0,0){$t_0$}}
\put(4,18.5){\makebox(0,0){$t_4$}}
\put(46.5,18.5){\makebox(0,0){$t_1$}}
\put(4,1.5){\makebox(0,0){$t_3$}}
\put(46.5,1.5){\makebox(0,0){$t_2$}}

\put(5,24){\makebox(0,0){$y_1$}}
\put(45,24){\makebox(0,0){$y_3$}}
\put(25,10){\makebox(0,0){$y_2$}}
\put(1,10){\makebox(0,0){$y_0$}}
\put(49,10){\makebox(0,0){$y_4$}}
\put(25,-3){\makebox(0,0){$z$}}

\end{picture} 
\end{center} 
\caption{Close-up of a brick %Figure~\ref{fig:brick wall}
} 
\label{fig:brick} 
\end{figure} 

For every vertex $t$ on the median of the wall,
the matrix $B(t)$ is obtained from $B$ by a simultaneous permutation of 
rows and columns. 
Now take a vertex on the outer boundary, 
say the vertex $t_2$ in Figure~\ref{fig:brick}. 
Then 
$$B(t_2) = \mu_3 (B(t_1)) = 
\left[\begin{array}{ccc} 
0 & -2 & 1 \\ 
2 & 0 & -1\\ 
-1 & 1 & 0 \\ 
\end{array}\right] \ .$$
Again, for every other vertex on the outer boundary, 
the corresponding matrix is obtained from $B(t_2)$ by a 
simultaneous permutation of rows and columns.
Substituting the matrices $B(t)$ into (\ref{eq:M_k(t)}), we generate 
all the exchange relations (\ref{eq:brick exchange}) 
(cf.~Figure~\ref{fig:brick}). 

To prove that our brick wall is indeed the exchange graph, 
it remains to show that all the cluster variables $y_m$, $w$ and $z$ 
are distinct;
recall that we view them as elements of the ambient field~$\FFcal$. 
Following the methodology of Theorem~\ref{th:denominators-rank 2}, it suffices
to show that all of them have different denominators in the Laurent expansion 
in terms of the initial cluster $\{y_1, y_2, y_3\}$.
We write the denominator of a cluster variable $y$ as 
$y_1^{d_1 (y)} y_2^{d_2 (y)} y_3^{d_3 (y)}$ 
(cf.~(\ref{eq:denominator})), and encode it by a vector 
$\delta (y) = (d_1 (y), d_2 (y), d_3 (y)) \in \ZZ^3$. 
In particular, we have
\begin{equation}
\label{eq:initial denominators}
\delta (y_1) = (-1,0,0), \,\, \delta (y_2) = (0,-1,0), \,\,
\delta (y_3) = (0,0,-1) \ .
\end{equation}
A moment's reflection shows that the vectors 
$\delta (y_m), \delta (w)$, and $\delta (z)$ 
satisfy the ``tropical version" of each 
of the exchange relations (\ref{eq:brick exchange}) 
(with multiplication replaced by vector addition, and addition replaced 
by component-wise maximum),
as long as the two vectors on the right-hand side are different.
For example, the first equation in (\ref{eq:brick exchange}) will take the form
$$
\delta (w) + \delta (y_{2m}) 
= \max(\delta (y_{2m-1}),\delta (y_{2m+1})) 
\ ,$$
provided $\delta (y_{2m-1}) \neq \delta (y_{2m+1})$;
here $\max$ stands for component-wise maximum. 
%There is an important caveat here: 
%the tropical version of each of the equations (\ref{eq:brick exchange}) holds 
%provided the two vectors in the right-hand side are different. 
All $\delta (y)$'s are uniquely determined 
from~(\ref{eq:initial denominators}) by recursive application of these conditions. 
As a result, we obtain, for every $m \geq 0$:
\[
\begin{array}{ccc}
&\delta (w) = (0,1,0) & \\[.1in]
\delta (y_{1-2m}) = (m-1, m, m) & &\delta (y_{2m+3}) = (m, m, m-1) \\[.1in]
\delta (y_{-2m}) = (m, m, m+1) & & \delta (y_{2m+4}) = (m+1, m, m) \\[.1in]
 & \delta (z) = (1,0,1)&
\end{array}
\]
(cf.~Figure~\ref{fig:brick wall exponents}, where each vector
$\delta (y)$ is shown within the corresponding region). 

\begin{figure}[ht] 
\setlength{\unitlength}{1.38pt}
\begin{center} 
\begin{picture}(260,63)(0,0)
%\thicklines
\multiput(0,10)(0,20){3}{\line(1,0){260}}
\multiput(10,30)(40,0){7}{\line(0,1){20}}
\multiput(30,10)(40,0){6}{\line(0,1){20}}

%\thinlines
\multiput(10,30)(40,0){7}{\circle*{2}}
\multiput(10,50)(40,0){7}{\circle*{2}}
\multiput(30,30)(40,0){6}{\circle*{2}}
\multiput(30,10)(40,0){6}{\circle*{2}}

\put(10,20){\makebox(0,0){$\begin{array}{c}y_{-4}\\(2,2,3)\end{array}$}}
\put(50,20){\makebox(0,0){$\begin{array}{c}y_{-2}\\(1,1,2)\end{array}$}}
\put(90,20){\makebox(0,0){$\begin{array}{c}y_0\\(0,0,1)\end{array}$}}
\put(130,20){\makebox(0,0){$\begin{array}{c}y_2\\(0,-1,0)\end{array}$}}
\put(170,20){\makebox(0,0){$\begin{array}{c}y_4\\(1,0,0)\end{array}$}}
\put(210,20){\makebox(0,0){$\begin{array}{c}y_6\\(2,1,1)\end{array}$}}
\put(250,20){\makebox(0,0){$\begin{array}{c}y_8\\(3,2,2)\end{array}$}}

\put(30,40){\makebox(0,0){$\begin{array}{c}y_{-3}\\(1,2,2)\end{array}$}}
\put(70,40){\makebox(0,0){$\begin{array}{c}y_{-1}\\(0,1,1)\end{array}$}}
\put(110,40){\makebox(0,0){$\begin{array}{c}y_1\\(-1,0,0)\end{array}$}}
\put(150,40){\makebox(0,0){$\begin{array}{c}y_3\\(0,0,-1)\end{array}$}}
\put(190,40){\makebox(0,0){$\begin{array}{c}y_5\\(1,1,0)\end{array}$}}
\put(230,40){\makebox(0,0){$\begin{array}{c}y_7\\(2,2,1)\end{array}$}}
%\put(190,40){\makebox(0,0){$\begin{array}{c}y_9\\(3,3,2)\end{array}$}}
%\put(230,40){\makebox(0,0){$\begin{array}{c}y_{11}\\(4,4,3)\end{array}$}}

%\put(10,20){\makebox(0,0){$y_0$}}
%\put(50,20){\makebox(0,0){$y_2$}}
%\put(90,20){\makebox(0,0){$y_4$}}
%\put(130,20){\makebox(0,0){$y_6$}}

%\put(30,40){\makebox(0,0){$y_1$}}
%\put(70,40){\makebox(0,0){$y_3$}}
%\put(110,40){\makebox(0,0){$y_5$}}

%\put(53,33){\makebox(0,0){$t_0$}}
%\put(73,33){\makebox(0,0){$t_1$}}
\put(130,0){\makebox(0,0){$\begin{array}{c}z\\(1,0,1)\end{array}$}}
\put(130,61){\makebox(0,0){$\begin{array}{c}w\\(0,1,0)\end{array}$}}

\end{picture} 
\end{center} 
\caption{The brick wall denominators} 
\label{fig:brick wall exponents} 
\end{figure}
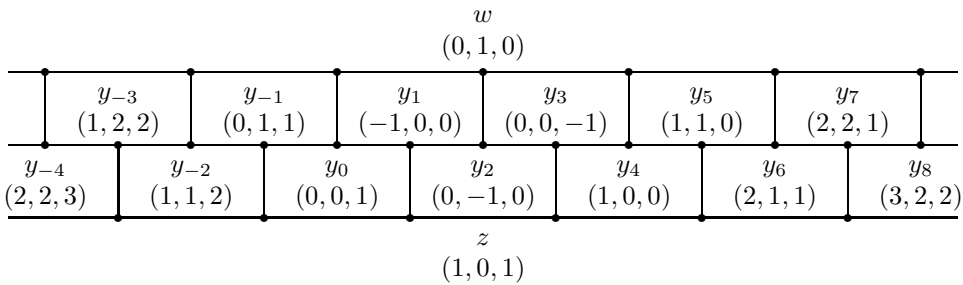 
}
\end{example}


\begin{thebibliography}{xxx} 
 
\bibitem{bfz96} 
A.~Berenstein, S.~Fomin and A.~Zelevinsky, 
Parametrizations of canonical bases and totally positive matrices, 
\textsl{Adv.\ Math.} {\bf 122} (1996), 49--149. 
 
 
\bibitem{bz93} 
A.~Berenstein and A.~Zelevinsky, 
String bases for quantum groups of type $A_r$, 
in I. M. Gelfand Seminar, 51--89, 
\textsl{Adv.\ Soviet Math.} {\bf 16}, Part 1, 
Amer.\ Math.\ Soc., Providence, RI, 1993. 
 
\bibitem{bz97} 
A.~Berenstein and A.~Zelevinsky, 
Total positivity in Schubert varieties, 
\textsl{Comment.\ Math.\ Helv.} {\bf 72} (1997), 128--166. 
 
\bibitem{bz01} 
A.~Berenstein and A.~Zelevinsky, 
Tensor product multiplicities, canonical bases and totally positive varieties, 
\textsl{Invent.\ Math.} {\bf 143} (2001), 77--128. 
 
%\bibitem{conway-guy} 
%J.~H.~Conway and R.~K.~Guy, \textsl{The book of numbers,} Copernicus, New York, 1996. 
 
 
%\bibitem{conway-coxeter} 
%J.~H.~Conway and H.~S.~M.~Coxeter, 
%Triangulated polygons and frieze patterns, 
%\textsl{Math.\ Gaz.} {\bf 57} (1973), no.~400, 87--94; 
%and \textsl{ibid.}, no.~401, 175--183. 
 
%\bibitem{fs-michigan} 
%S.~Fomin and M.~Shapiro, 
%Stratified spaces formed by totally positive varieties, 
%\textsl{Michigan Math.\ J.} {\bf 48}  (2000), 253--270. 
 
\bibitem{fz-jams} 
S.~Fomin and A.~Zelevinsky, 
Double Bruhat cells and total positivity, 
\textsl{J.\ Amer.\ Math.\ Soc.} {\bf 12} (1999), 335--380. 
 
 
\bibitem{fz-intelligencer} 
S.~Fomin and A.~Zelevinsky, 
Total positivity: tests and parametrizations, 
\textsl{Math.\ Intelligencer} {\bf 22} (2000), no.~1, 23--33. 
 
\bibitem{fz-proceedings} 
S.~Fomin and A.~Zelevinsky, 
Totally nonnegative and oscillatory elements in semisimple groups, 
\textsl{Proc.\ Amer.\ Math.\ Soc.} {\bf 128} (2000), 3749--3759. 
 
\bibitem{fz-Laurent} 
S.~Fomin and A.~Zelevinsky, 
The Laurent phenomenon, %for cluster exchange recurrences, 
in preparation. 
 
 
%\bibitem{gale-intelligencer} 
%D.~Gale, Mathematical Entertainments: The strange and surprising saga 
%of the Somos sequences, \emph{Math.\ Intelligencer} {\bf 13} (1991), 
%no.~1, 40--43. 
 
\bibitem{gz86} 
I.~M.~Gelfand and A.~Zelevinsky, 
Canonical basis in irreducible representations of $gl_3$ and its applications, 
in: \textsl{Group theoretical methods in physics, Vol.~II} (Jurmala, 1985), 
127--146, VNU Sci.\ Press, Utrecht, 1986. 
 
%\bibitem{guy-unsolved} 
%R.~K.~Guy, \textsl{Unsolved problems in number theory,} 2nd edition, 
%Springer-Verlag, New York, 1994. 
 
\bibitem{kac} 
V.~Kac, Infinite dimensional Lie algebras, 3rd edition, 
Cambridge University Press, 1990. 
 
\bibitem{kungrota84} 
J.~Kung and G.-C.~Rota, The invariant theory of binary forms, 
\textsl{Bull.\ Amer.\ Math.\ Soc.\ (N.S.)} {\bf 10} (1984), 27--85. 
 
\bibitem{lekzel98} 
B.~Leclerc and A.~Zelevinsky, 
Quasicommuting families of quantum Plucker coordinates, 
\textsl{Amer.\ Math.\ Soc.\ Transl.\ (2),} Vol.~181, 
Kirillov's Seminar on Representation Theory, 85--108, 
Amer.\ Math.\ Soc., Providence, RI, 1998. 
 
 
\bibitem{lu1} G.~Lusztig, 
Canonical bases arising from quantized enveloping algebras, 
\textsl{J.~Amer.\ Math.\ Soc.} {\bf 3} (1990), 447--498. 
 
\bibitem{lusztig-quantum} 
G.~Lusztig, {\sl Introduction to quantum groups}, 
Progress in Mathematics {\bf 110}, 
Birkh\"auser, 1993. 
 
\bibitem{lusztig-reductive} 
G.~Lusztig, Total positivity in reductive groups, 
in: {\sl Lie theory and geometry: 
in honor of Bertram Kostant, Progress in Mathematics} {\bf 123}, 
Birkh\"auser, 1994. 
 
 
 
\bibitem{rz} 
V.~Retakh and A.~Zelevinsky, 
The base affine space and canonical bases in irreducible representations of the group $Sp(4)$, 
\textsl{Soviet Math.\ Dokl.} {\bf 37} (1988), no.~3, 618--622. 
 
%\bibitem{robinson} 
%R.~M.~Robinson, Periodicity of Somos sequences, 
%\textsl{Proc.\ Amer.\ Math.\ Soc.} {\bf 116} (1992), 613--619. 
 
\bibitem{ssvz} 
B.~Shapiro, M.~Shapiro, A.~Vainshtein and A.~Zelevinsky, 
Simply-laced Coxeter groups and groups generated by symplectic transvections, 
\textsl{Michigan Math.\ J.} {\bf 48}  (2000), 531--551. 
 
 
\bibitem{stur} 
B.~Sturmfels, 
Algorithms in invariant theory. Texts and Monographs in Symbolic Computation. 
Springer-Verlag, Vienna, 1993. 
 
\bibitem{zamolodchikov} 
Al.~B.~Zamolodchikov, 
On the thermodynamic Bethe ansatz equations for reflectionless $ADE$ 
scattering theories, 
\textsl{Phys.\ Lett.~B} \textbf{253} (1991), 391--394. 
 
 
\bibitem{z-imrn} 
A.~Zelevinsky, 
Connected components of real double Bruhat cells, 
\textsl{Intern.\ Math.\ Res.\ Notices} 2000, No. 21, 1131--1153. 
 
\end{thebibliography}
\end{document}